\documentclass[11pt,a4paper,reqno]{amsart}
\usepackage{stmaryrd}
\usepackage[english]{babel}
\usepackage{pst-grad} 
\usepackage{pst-plot} 
\usepackage{pstricks}
\usepackage{amsmath,amssymb,mathrsfs,amsthm, mathtools}
\usepackage{tikz} 
\usepackage{mathcomp,wasysym}  
\usepackage{graphicx,subfigure}  
\usepackage[all]{xy} \xyoption{arc} \xyoption{color}
\usepackage{epsfig}
\usepackage{cite}
\usepackage[a4paper,
left=2.5cm, right=2.5cm,
top=3cm, bottom=3cm]{geometry}

\newcommand{\T}{\mathbb T}

\normalsize
\normalsize
\newcommand{\norm}[1]{\left\lVert#1\right\rVert} 
\DeclarePairedDelimiter\abs{\lvert}{\rvert}%

\usepackage{hyperref}
\hypersetup{
colorlinks=true,
linktoc=all,
citecolor=black,
linkcolor=blue,
}
\newtheorem{satz}{Proposition}[section]
\newtheorem{lem}[satz]{Lemma} 
\newtheorem{theorem}{Theorem}
\newtheorem*{theorem*}{Theorem}

\newtheorem{lemma}[satz]{Lemma} 

\theoremstyle{definition}

\newtheorem{definition}{Definition}

\title[Traveling wave solutions for cold plasmas]{On the existence of traveling wave solutions for cold plasmas}

\author[D. Alonso-Or\'{a}n]{Diego Alonso-Or\'{a}n}
\email{dalonsoo@ull.edu.es}
\address{Departamento de An\'{a}lisis Matem\'{a}tico y Instituto de Matem\'{a}ticas y Aplicaciones (IMAULL), Universidad de La Laguna C/. Astrof\'{i}sico Francisco S\'{a}nchez s/n, 38200 - La Laguna, Spain.}
\author[A. Dur\'{a}n]{Angel Dur\'{a}n}
\email{angeldm@uva.es}
\address{Applied Mathematics Department, University of Valladolid, 47011 Valladolid, Spain}
\author[R. Granero-Belinch\'{o}n]{Rafael Granero-Belinch\'{o}n}
\email{rafael.granero@unican.es}
\address{Departamento  de  Matem\'aticas,  Estad\'istica  y  Computaci\'on,  Universidad  de Cantabria.  Avda.  Los  Castros  s/n,  Santander,  Spain.}

\begin{document}
\begin{abstract}
The present paper is concerned with the existence of traveling wave solutions of the asymptotic model, derived by the authors in a previous work, to approximate the unidirectional evolution of a collision-free plasma in a magnetic field. First, using bifurcation theory, we can rigorously prove the existence of periodic traveling waves of small amplitude. Furthermore, our analysis also evidences the existence of different type of traveling waves. To this end, we present a second approach based on the analysis of the differential system satisfied by the traveling-wave profiles, the existence of equilibria, and the identification of associated homo-clinic and periodic orbits around them. The study makes use of linearization techniques and numerical computations to show the existence of different types of traveling-wave solutions, with monotone and non-monotone behaviour and different regularity, as well as periodic traveling waves. 
\end{abstract}



\maketitle
{\small
\tableofcontents}

\allowdisplaybreaks

\section{Introduction}
The motion of a magnetized cold plasma consisting of singly-charged particles can be described by the following hyperbolic–hyperbolic–elliptic system of PDEs \cite{berezin1964theory,gardner1960similarity}
\begin{subequations}\label{eq:1}
\begin{align}
n_t+(un)_x&=0,\\
u_t+uu_x+\frac{bb_x}{n}&=0,\\
b-n-\left(\frac{b_x}{n}\right)_x&=0,
\end{align}
\end{subequations}
where $n$ represents the number density of ions, $u$ is the ion velocity and $b$ the magnetic field. In \eqref{eq:1} the unknowns are functions of $(x,t)\in\mathbb{R}\times [0,\infty)$.  \medskip

The system \eqref{eq:1} was introduced in \cite{gardner1960similarity} to study hydromagnetic waves traveling across a magnetic field. Later studies regarding the  oblique propagation of hydromagnetic waves were investigated in \cite{berezin1964theory, kakutani1968reductive}. The rigorous justification of the KdV limit of \eqref{eq:1} is provided in \cite{pu2019kdv}. The first and third author showed in \cite{ag2024} the local well-posedness of classical solutions to \eqref{eq:1} with initial data in $(\rho_{0}-1,u_{0})\in H^{2}(\mathbb{R})\times H^{3}(\mathbb{R})$. Later, in \cite{baechoikwon2024}, Bae, Choi and Kwon demonstrate that solutions to \eqref{eq:1} blow-up in finite time for a certain class of initial data.

\medskip Recently in \cite{adg2024}, by means a multi-scale expansion approach, the authors derived asymptotic models of \eqref{eq:1}. The first one is the nonlocal Boussinesq system
\begin{subequations}\label{eq:Boussinesq}
\begin{align}
h_t+(hv)_x+v_x&=0,\\
v_t+vv_x+\left[\mathscr{L},\mathscr{N}h\right]h+\mathscr{N}h&=0,
\end{align}
\end{subequations}
where the nonlocal operators are given by
\begin{equation}\label{operators:nonlocal}
\mathscr{L}=-\partial_x^2(1-\partial_x^2)^{-1},\;
\mathscr{N}=\partial_x(1-\partial_x^2)^{-1}\;(\text{so }\partial_x\mathscr{N}=-\mathscr{L}),
\end{equation}
with Fourier symbols
$$
\widehat{\mathscr{L}h}(\xi)=\frac{\xi^2}{1+\xi^2}\hat{h}(\xi),
\; \widehat{\mathscr{N}h}(\xi)=\frac{i\xi}{1+\xi^2}\hat{h}(\xi),
$$
and where $\left[\mathscr{L},\cdot\right]\cdot$ denotes the commutator
\begin{equation*}\label{adg4b}
\left[\mathscr{L},f\right]g=\mathscr{L}(fg)-f\mathscr{L}g.
\end{equation*}
In (\ref{eq:Boussinesq}) $h$ and $v$ represent, respectively, second-order approximations in a multi-scale expansion of the ionic density and ionic velocity variables. An additional assumption concerning the first-order terms of the expansions leads to the single, unidirectional, nonlocal wave equation for $h$ given by
\begin{equation}\label{wave:uni}
h_{t}=-\frac{1}{2}\left(3 hh_{x}-[\mathscr{L},\mathscr{N}h]h - \mathscr{N}h-h_{x} \right).
\end{equation}
Furthermore, besides the derivation of \eqref{eq:Boussinesq} and \eqref{wave:uni}, the well-posedness of the corresponding initial-value problems (ivp's) and additional properties, such as the Hamiltonian structure, existence of conserved quantities and the formation of singularities in finite time are also studied in \cite{adg2024}. \medskip

Our next step in the analysis of (\ref{eq:Boussinesq}) and (\ref{wave:uni}) is concerned with the existence of traveling wave solutions. These are solutions of permanent form and traveling with some constant velocity. Their determination transforms (\ref{eq:Boussinesq}) and (\ref{wave:uni}) into the corresponding ordinary differential equations for the wave profiles, depending on the speed. If the profiles are known to go to zero asymptotically (surely along with their derivatives up to some order), that is, if the profiles are localized, they are usually called solitary waves. The present paper is focused on the unidirectional equation (\ref{wave:uni}), while the system (\ref{eq:Boussinesq}) will be considered in a forthcoming work. \medskip

The study of the existence of solitary wave solutions, specially for models in water wave theory, may be accomplished with several, classical theories, e.~g. \cite{BenjaminBB1990,IK,Lions,Lombardi3,Toland1986}. In the first result provided in this manuscript we are able to rigorously show the existence of periodic traveling waves of small amplitude and $O(1)$ velocity, cf. \cite{adg2024}. The proof hinges on non-linear bifurcation techniques via Crandall-Rabinowitz theorem. The bifurcation approach to prove the existence of traveling waves for different equations has been used by different authors \cite{AGG,Constantin-Varvaruca11,Lenells1,Lenells2}. \medskip

Furthermore, we will also make some formal arguments and simulations that evidence the existence of traveling waves of different type. More precisely, we show a picture on the existence of different types of traveling-wave solutions of  (\ref{wave:uni}) (including, but not only, those of solitary type). This approach is based on the following contributions:
\begin{itemize}
\item The identification of a reformulation of (\ref{wave:uni}) with only local terms.
\item The corresponding differential (ode) systems for the traveling wave profiles are shown to have several remarkable properties which allow to determine and classify different families of equilibria, according to the speed of the wave, used as bifurcation parameter.
\item In addition, the reversible character of the ode systems can be used to analyze the dynamics of some of these equilibria from the literature on the emergence and homo-clinic and periodic orbits via linearization and Normal Form theory for reversible vector fields, \cite{Champ,IK,HaragusI}. Their application reveals the existence of different types of traveling wave solutions, with monotone and non-monotone behaviour and different regularity, as well as smooth solitary waves and periodic traveling waves.
\item The discussion is supported by some experiments from an efficient numerical procedure to compute approximations to the profiles. The numerical results suggest some additional properties of the waves, concerning amplitude, speed and decay at infinity.
\end{itemize}
\medskip

\subsection*{Plan of the paper} The paper is structured as follows. In section \ref{notation:aux} we fix the notation and present some auxiliary results on bifurcation theory. In section \ref{sec:teo} the existence of smooth periodic traveling waves of small amplitude is proved via a bifurcation argument.  Section \ref{sec2} is devoted to derive the local formulation of (\ref{wave:uni}) and to identify the families of equilibria of the corresponding ode system satisfied by the traveling waves. For two of these equilibria and from the reversibility property of the differential equations, the existence of traveling wave solutions as homo-clinic and periodic orbits around them is discussed in section \ref{sec3}. The discussion is based on the application of normal form theory of the corresponding reversible vector fields, \cite{HaragusI}, with the help of literature for similar equations, (cf., among others, \cite{Champ,ChampS,ChampT,DDS1,IK} and references therein), and the numerical experiments of generation of approximate waves. Some concluding remarks are outlined in section \ref{sec4}. The numerical procedure to compute approximations of traveling wave profiles is described in appendix \ref{appA}.

\subsection*{Acknowledgments}
The authors gratefully acknowledge Claudia Garc\'ia for fruitful discussions regarding bifurcation theory and Daniel S\'anchez Sim\'on del Pino for pointing out the reference of Lemma 2.1. D.A-O is supported by the fellowship of the Santander-ULL program. D.A-O and R. G-B are also supported by the project “An\'alisis Matem\'atico Aplicado y Ecuaciones Diferenciales” Grant PID2022- 141187NB-I00 funded by MCIN/ AEI and acronym “AMAED”. This paper is part of the project PID2022-141187NB-I00 with acronym \emph{AMAED} funded by MICIU/AEI /10.13039/501100011033 and by FEDER, UE/AEI /10.13039/501100011033 / FEDER, UE. A.D. is supported by the Spanish Agencia Estatal de Investigación under Research Grant PID2023-147073NB-100.

\section{Notation and auxiliary results}\label{notation:aux}
In this section, we fix the notation used throughout the paper and recall some classical results on bifurcation theory. For a periodic function $f$ with values in $\mathbb{R}$ we define the H\"older norms as 
\begin{align*}
\norm{f}_{C^0(\T)}&=\sup_{x\in \mathbb{\T}^1} \abs{f}, \; \norm{f}_{C^k(\T)}=\norm{f}_{C^0(\T)}+\sum_{\ell=1}^k \norm{\partial_x^\ell f}_{C^0(\T)},\; k\in \mathbb{N},\\
\norm{f}_{C^{0,\alpha}(\T)}&=\norm{f}_{C^0(\T)}+\sup_{x_1,x_2\in (\T)} \frac{\abs{f(x_1)-f(x_2)}}{\abs{x_1-x_2}^{\alpha}}, \; 0<\alpha<1,\\
\norm{f}_{C^{k,\alpha}(\T)}&=\norm{f}_{C^{k-1}(\T)}+\norm{\partial_x^k f}_{C^\alpha(\T)}, \; k\in \mathbb{N},\; 0<\alpha<1.
\end{align*}
The Banach space of continuous functions for which the above norms are finite will be denoted $C^{k}(\T)$ and $C^{k,\alpha}(\T)$. \medskip

We denote by $\mathscr{Q}:=(\mathsf{I}-\partial_{xx}^{2})^{-1}$ the Helmholtz operator with Fourier symbol
\[ \widehat{\mathscr{Q}f}(\xi)=\frac{1}{1+\xi^{2}}\widehat{f}(\xi).\]
Acting on square integrable periodic functions $f$, the Helmholtz operator has the representation
\begin{equation}\label{Helmholtz}
 \mathscr{Q}f(x)=[(\mathsf{I}-\partial_{xx}^{2})^{-1}]f(x)=[\mathsf{G}_{\T}\star f](x),
 \end{equation}
where the Green function $\mathsf{G}_{\T}$ is explicitly given by
\begin{equation}\label{Helmholtz:G}
\mathsf{G}_{\T}(x)=\frac{\cosh(x-2\pi [\frac{x}{2\pi}]-\pi)}{2\sinh(\pi)}.
\end{equation}\medskip
We have the following bound whose proof can be found in \cite[Theorem 4, \S 4.4]{SteinDiff}:
\begin{lem}\label{lemma:bound}
The operator $\mathscr{Q}$ maps $C^{k,\alpha}(\mathbb{T})$ isomorphically onto $C^{k+2,\alpha}(\mathbb{T})$. More precisely, there exists a constant $C>0$ such that
\begin{equation}\label{bound:helmholtz}
\norm{\mathscr{Q}f}_{C^{k+2,\alpha}(\mathbb{T})}\leq C \norm{f}_{C^{k,\alpha}(\mathbb{T})}, \quad f\in C^{k,\alpha} (\mathbb{T}), \quad k\in \mathbb{N}\cup\{0\}.
\end{equation}
\end{lem}
Therefore as a direct consequence the non-local operators in \eqref{operators:nonlocal} written as
\begin{equation}\label{relation:nonlocal:Helm}
 \mathscr{L}=-\partial_{x}^{2}\mathscr{Q}, \; \mathscr{N}=\partial_x \mathscr{Q}.
 \end{equation} 
enjoyed the following estimates
\begin{equation}\label{bound:nonloca}
\norm{\mathscr{L} f}_{C^{k,\alpha}(\mathbb{T})}\leq C\norm{f}_{C^{k,\alpha}(\mathbb{T})}, \;  \norm{\mathscr{N} f}_{C^{k+1,\alpha}(\mathbb{T})}\leq C\norm{f}_{C^{k,\alpha}(\mathbb{T})}.
\end{equation}
The previous bounds \eqref{bound:helmholtz}-\eqref{bound:nonloca} are particular cases of a more general theory of mapping properties for pseudo-differential operators in Besov spaces, cf. \cite[Theorem 6.19, \S 6.6]{Abels}. \medskip

Next, let us revisit the Crandall-Rabinowitz Theorem, an essential tool in bifurcation theory that we will use to show the existence of smooth periodic traveling waves.  Before we dive in, let us first go over the following definition:

\begin{definition}[Fredholm operator]
Let $X$ and $Y$ be  two Banach spaces. A continuous linear mapping $T:X\rightarrow Y,$  is a  Fredholm operator if it fulfills the following properties,
\begin{enumerate}
\item $\textnormal{dim Ker}\,  T<\infty$,
\item $\textnormal{Im}\, T$ is closed in $Y$,
\item $\textnormal{codim Im}\,  T<\infty$.
\end{enumerate}
\end{definition}
The integer $\textnormal{dim Ker}\, T-\textnormal{codim Im}\, T$ is called the Fredholm index of $T$. Moreover, we also remark that the index of a Fredholm operator remains unchanged under compact perturbations,  cf. \cite{Kato, Kiel}. Next, let us state the classical Crandall-Rabinowitz theorem \cite{Crandall-Rabi} which reads as follows:
\begin{theorem}[Crandall-Rabinowitz Theorem]\label{CR}
    Let $X, Y$ be two Banach spaces, $V$ be a neighborhood of $0$ in $X$ and $F:\mathbb{R}\times V\rightarrow Y$ be a function with the properties,
    \begin{enumerate}
        \item $F(\lambda,0)=0$ for all $\lambda\in\mathbb{R}$.
        \item The partial derivatives  $\partial_\lambda F$, $\partial_fF$ and  $\partial_{\lambda}\partial_fF$ exist and are continuous.
        \item The operator $\partial_f F(\lambda_{0},0)$ is Fredholm of zero index and $\textnormal{Ker}(\partial_f F(\lambda_{0},0))=\langle f_0\rangle$ is one-dimensional. 
                \item  Transversality assumption: $\partial_{\lambda}\partial_fF(\lambda_{0},0)f_0 \notin \textnormal{Im}(\partial_fF(\lambda_{0},0))$.
    \end{enumerate}
    If $Z$ is any complement of  $\textnormal{Ker}(\partial_fF(\lambda_{0},0))$ in $X$, then there is a neighborhood  $U$ of $(\lambda_{0},0)$ in $\mathbb{R}\times X$, an interval  $(-a,a)$, and two continuous functions $\Phi:(-a,a)\rightarrow\mathbb{R}$, $\beta:(-a,a)\rightarrow Z$ such that $\Phi(0)=\lambda_{0}$ and $\beta(0)=0$ and
    $$F^{-1}(0)\cap U=\{(\Phi(s), s f_0+s\beta(s)) : |s|<a\}\cup\{(\lambda,0): (\lambda,0)\in U\}.$$
\end{theorem}
In this context, we will say that $\lambda_{0}$ is an eigenvalue of $F$.

\section{Existence of periodic traveling wave solutions}\label{sec:teo}

This section is devoted to show the existence of smooth periodic traveling waves of small amplitude. The precise statement of result reads:
\begin{theorem}\label{main:theorem}
For any $m\geq 1$ there exists a one dimensional curve $s\mapsto (c_s, \varphi_s)$, with $s\in I$, such that
$$
h(x)=\varphi_s(x)\in C^{1,\alpha}([0,2\pi],\mathbb{R}),
$$
is a $m-$fold traveling wave solution to \eqref{wave:uni} with constant speed $c_s$.
\end{theorem}

\subsection{Formulation and functional spaces} We look for periodic traveling waves for $h$ and hence we make the Ansatz 
$$
h(x,t)=\varphi(x-ct),
$$
for some speed  $c\in\mathbb{R}$. Hence, plugging it into \eqref{wave:uni} we find the problem
\begin{equation}\label{wave:unitraveling}
-c\varphi'=-\frac{1}{2}\left(3 \varphi\varphi'-[\mathscr{L},\mathscr{N}\varphi]\varphi - \mathscr{N}\varphi-\varphi' \right).
\end{equation}
As a consequence, the equation reduces to
$$
F[c,\varphi](\xi)=0, \quad \xi\in[-\pi,\pi],
$$
where
\begin{align}
F[c,\varphi](\xi)&=\frac{1}{2}\left(3 \varphi(\xi)\varphi'(\xi)-[\mathscr{L},\mathscr{N}\varphi(\xi)]\varphi(\xi) - \mathscr{N}\varphi(\xi)-\varphi'(\xi)\right)-c\varphi'(\xi). \label{defiF}
\end{align}
We observe that, regardless of the value of $c$ we have the following line of trivial solutions
$$
F[c,0]=0.
$$
Following \cite{AGG}, we define the functional spaces
\begin{align*}
X:=\left\{h\in C^{1,\alpha}([0,2\pi],\mathbb{R}),\quad h(\xi)=\sum_{k\geq 1}h_k\cos(k\xi)\text{ with norm }\|h\|_{X}=\|h\|_{C^{1,\alpha}}\right\},\\
Y:=\left\{h\in C^{0,\alpha}([0,2\pi],\mathbb{R}),\quad h(\xi)=\sum_{k\geq 1}h_k\sin(k\xi)\text{ with norm }\|h\|_{Y}=\|h\|_{C^{0,\alpha}}\right\}.
\end{align*}
Furthermore, we also introduce the space 
\[ Z:=\left\{h\in C^{2,\alpha}([0,2\pi],\mathbb{R}),\quad h(\xi)=\sum_{k\geq 1}h_k\sin(k\xi)\text{ with norm }\|h\|_{Z}=\|f\|_{C^{2,\alpha}}\right\}. \]
It is straightforward to check that the embedding  $Z\hookrightarrow Y$ is compact. \medskip

\subsection{Spectral and tranversality properties of the linearized operator} \label{subsec:spectral}
In this subsection, we will check that the hypothesis of the Crandall-Rabinowitz Theorem \ref{CR} are satisfied. To start with, let us show that the operator $F:\mathbb{R}\times X\rightarrow Y$ given in \eqref{defiF} is well-defined and $\mathscr{C}^1(\mathbb{R}\times X\rightarrow Y)$. We first observe that if $\varphi$ is an even function, $\varphi(\xi)=\varphi(-\xi)$, then $F[c,\varphi]$ is an odd function, i.e.,
\[ F[c,\varphi](\xi)=-F[c,\varphi](-\xi).\]
To that purpose, we first notice that $\varphi'(-\xi)=-\varphi'(\xi)$. Moreover, using the representation \eqref{Helmholtz}-\eqref{Helmholtz:G} for the Helmholtz operator, we easily check that
\[ \mathscr{Q}\varphi(\xi)=\mathscr{Q}\varphi(-\xi).\] 
Invoking \eqref{relation:nonlocal:Helm} we find that
\[ \mathscr{N}\varphi(\xi)=\mathscr{Q}\varphi'(\xi)=-\mathscr{N}\varphi(-\xi), \quad \mathscr{L}\varphi(\xi)=(\mathsf{I}-\mathscr{Q})\varphi(\xi)=\mathscr{L}\varphi(-\xi), \]
and hence $[\mathscr{L},\mathscr{N}\varphi]\varphi(\xi)=-[\mathscr{L},\mathscr{N}\varphi]\varphi(-\xi).$ Thus, by combining the previous identities and recalling \eqref{defiF} we conclude that $F[c,\varphi](\xi)=-F[c,\varphi](-\xi)$ satisfying the symmetry property. \medskip

In the following let us show that $F:\mathbb{R}\times X\rightarrow Y$ is well-defined. It is straightforward to check that
\[ \norm{-\frac{1}{2}\left(\mathscr{N}\varphi+\varphi'\right)-c\varphi'}_{Y}\leq C\norm{\varphi}_{X}. \]
Similarly, the non-linear terms can be bounded by
\[ \norm{\frac{1}{2}\left( 3\varphi\varphi'-[\mathscr{L},\mathscr{N}\varphi]\varphi\right)}_{Y}\leq C\norm{\varphi}_{X}\norm{\varphi}_{Y}+\norm{\varphi}_{Y}^{2}, \]
where we have used estimate \eqref{bound:nonloca} and the fact that 
\[ \norm{fg}_{Y}\leq \norm{f}_{Y}\norm{g}_{Y}.\]
Thus, since $X\hookrightarrow Y$, we conclude that
 $$
\|F[c,\varphi]\|_{Y}\leq C\left(\|\varphi\|_{X}^2+\|\varphi\|_{X}\right).
$$
In order to prove that $F\in \mathscr{C}^1(\mathbb{R}\times X\rightarrow Y)$ it is sufficient to show that the Gateaux derivative of $F$ verifies
\begin{equation}\label{FC1}
\|\partial_\varphi F[c,\varphi_1]f-\partial_\varphi F[c,\varphi_2]f\|_{Y}\leq C\|f\|_{X}\|\varphi_1-\varphi_2\|_{X}. 
\end{equation}
We compute the Gateaux derivative of $F$ and find that
$$
\partial_\varphi F[c,\varphi]f=\frac{3}{2}\varphi'f+\frac{3}{2}\varphi f'-[\mathscr{Q},\mathscr{N}f]\varphi-[\mathscr{Q},\mathscr{N}\varphi]f-\frac{1}{2}\mathscr{N}f-\frac{1}{2}f'-cf'.
$$
Repeating the same estimates by making use of estimates \eqref{bound:nonloca} and Lemma \ref{lemma:bound} to deal with the Helmholtz operator $\mathscr{Q}$ we can easily check that \eqref{FC1} is indeed satisfied. Hence, we can conclude that the Gateaux derivative is continuous (indeed, it is Lipschitz) and then we can ensure the existence and continuity of the Fr\'echet derivative.\medskip

Next, we analyze the linearized operator at the trivial solution which is given by
\begin{equation*}
\partial_\varphi F[c,0]f(\xi)=-\frac{1}{2}\mathscr{N}f(\xi)-\frac{1}{2}f'(\xi)-cf'(\xi):=\mathsf{L}(f)(\xi)+\mathsf{K}(f)(\xi),
\end{equation*}
where
\[ \mathsf{L}(f)(\xi)=-\frac{1}{2}f'-cf, \ \mbox{ and }  \mathsf{K}(f)(\xi)=-\frac{1}{2}\mathscr{N}f.\]
The principal part of the operator $\mathsf{L}f: X\rightarrow Y $ is an isomorphism and hence has zero index. Moreover, using estimate \eqref{bound:nonloca}, we infer that the operator $\mathsf{K}f:X\rightarrow Z$ is continuous and the embedding $Z\hookrightarrow Y$ is compact. Therefore, since the index of a Fredholm operator remains unchanged due to compact perturbations we conclude that $\partial_\varphi F[c,0]f$ is a Fredholm operator of zero index. \medskip

To conclude we characterize the kernel of the linear operator. For $f\in X$ we find that
$$
\partial_\varphi F[c,0]f=\sum_{k=1}^\infty f_k\sin(kx)\left(-\frac{k}{2(1+k^2)}-\frac{k}{2}-ck\right).
$$
Thus, for 
$$
c_k=-\frac{1}{2(1+k^2)}-\frac{1}{2}
$$
we have that
$$
Ker(\partial_\varphi F[c,0])=\text{span}(\cos(kx))
$$
and, recalling that the linearized operator is Fredholm operator of zero index,
$$
Y/Img(\partial_\varphi F[c,0])=\text{span}(\sin(kx)).
$$
The transversality condition is then satisfied because
$$
\partial_c \partial_\varphi F[c_k,0]f=\sum_{k=1}^{\infty} f_k\sin(kx)\left(-k\right),
$$
and, if $f\in Ker(\partial_\varphi F[c,0])$, then
$$
\partial_c \partial_\varphi F[c_k,0]h=-kf_k\sin(kx)\notin Img(\partial_\varphi F[c,0]).
$$
\subsection{Proof of Theorem \ref{main:theorem}}
We have checked that the Crandall-Rabinowitz theorem can be applied in our equation \eqref{wave:unitraveling}. Fix $m\geq 1$. In order to prove Theorem \ref{main:theorem}, let us introduce the symmetry $m$ in the spaces. For that, let us define
\begin{align*}
X_m:=\left\{f\in C^{1,\alpha}([0,2\pi]),\quad f(\xi)=\sum_{k\geq 1}f_k\cos(mk\xi)\text{ with norm }\|f\|_{X_m}=\|f\|_{C^{1,\alpha}}\right\},\\
Y_m:=\left\{f\in C^{0,\alpha}([0,2\pi]),\quad f(\xi)=\sum_{k\geq 1}f_k\sin(mk\xi)\text{ with norm }\|f\|_{Y_m}=\|f\|_{C^{0,\alpha}}\right\},
\end{align*}
for any $m\geq 1$. The fact that the operator $F:\mathbb{R}\times X_{m}\to Y_{m}$ is well-defined and $\mathscr{C}^1(\mathbb{R}\times X_{m}\rightarrow Y_{m})$ can be checked by repeating the computations of Subsection \ref{subsec:spectral}. However, we have to show that the $m$-fold symmetry property holds. More precisely, we have to proof that if $\varphi(\xi+\frac{2\pi}{m})=\varphi(\xi)$ then
\[ F[c,\varphi](\xi+\frac{2\pi}{m})=F[c,\varphi](\xi).\]
Note that if $\varphi$ has $m$-fold symmetry, then all the derivatives of $\varphi$ also enjoy the same symmetry property. Moreover, recalling that $\mathscr{Q}$ is a convolution operator \eqref{Helmholtz}, the symmetry is also satisfied and thus $F[c,\varphi](\xi+\frac{2\pi}{m})=F[c,\varphi](\xi)$ follows. The rest of the arguments can be argued similarly as in the previous Subsection \ref{subsec:spectral}. Thus, Crandall-Rabinowitz theorem can be applied obtaining Theorem \ref{main:theorem}.

\section{Traveling wave solutions of the nonlocal wave equation}
\label{sec2}
In this section the existence of traveling wave solutions of the unidirectional model is analyzed. Note first, \cite{adg2024}, that
\eqref{wave:uni} can be written in the alternative form
\begin{equation}\label{adg8}
h_{t}+\partial_{x}\left(\frac{3}{4}h^{2}+\frac{1}{2}\mathscr{N}(h\mathscr{N}h)-\frac{1}{4}(\mathscr{N}h)^{2}-\frac{1}{2}\mathscr{Q}h-\frac{h}{2}\right)=0.
\end{equation}
An additional property that will be used below is a formulation of \eqref{wave:uni} involving local terms: if
\begin{eqnarray}
\mathscr{J}=\mathscr{Q}^{-1}, h=\mathscr{J}u,\label{loc}
\end{eqnarray} 
then \eqref{adg8} can be written as
\begin{equation}\label{adg8b}
\mathscr{J}u_{t}+\partial_{x}\left(\frac{3}{4}\left(\mathscr{J}u\right)^{2}+\frac{1}{2}\partial_{x}\mathscr{J}^{-1}\left((\mathscr{J}u)\partial_{x}u\right)
-\frac{1}{4}(\partial_{x}u)^{2}-\frac{1}{2}u-\frac{\mathscr{J}u}{2}\right)=0.
\end{equation}
For traveling wave solutions $h=h(x-c_{s}t)$ of (\ref{adg8}), the profiles $h=h(X), X=x-c_{s}t$ must satisfy
\begin{equation}\label{adg9}
-c_{s}h+\left(\frac{3}{4}h^{2}+\frac{1}{2}\mathscr{N}(h\mathscr{N}h)-\frac{1}{4}(\mathscr{N}h)^{2}-\frac{1}{2}\mathscr{Q}h-\frac{h}{2}\right)+g=0,
\end{equation}
with $g$ constant. In terms of the variable $u$, defined in (\ref{loc}), (\ref{adg9}) yields
\begin{equation}\label{adg9b}
-c_{s}\mathscr{J}^{2}u+\left(\frac{3}{4}\mathscr{J}\left(\mathscr{J}u\right)^{2}+\frac{1}{2}\partial_{x}\left((\mathscr{J}u)\partial_{x}u\right)
-\frac{1}{4}\mathscr{J}(\partial_{x}u)^{2}-\frac{1}{2}\mathscr{J}u-\frac{\mathscr{J}^{2}u}{2}\right)=0,
\end{equation}
which can be written as
\begin{eqnarray}
\left(\frac{3}{2}\mathscr{J}u-\widetilde{c}_{s}\right)u''''-\left(2\widetilde{c}_{s}+\frac{1}{2}\right)u''+\left(\widetilde{c}_{s}+\frac{1}{2}\right)u&&\nonumber\\
=\widetilde{f}(u,u',u'',u''')=f(u,u',u'',u''')+g&&\nonumber\\
=\frac{1}{2}\left(\frac{3}{2}u^{2}+\frac{7}{2}(u')^{2}+uu''-\frac{3}{2}(u'')^{2}-6u'u'''+3(u''')^{2}\right)+g,&&\label{adg9c}
\end{eqnarray}
where $\widetilde{c}_{s}=c_{s}+\frac{1}{2}$.
Equation (\ref{adg9b}) can also be formulated as a first-order system for $Y=(y_{1},y_{2},y_{3},y_{4})^{T}=(u,u',u'',u''')^{T}$ as
\begin{eqnarray}
\frac{dy_{1}}{dX}&=&y_{2},\nonumber\\
\frac{dy_{2}}{dX}&=&y_{3},\nonumber\\
\frac{dy_{3}}{dX}&=&y_{4},\nonumber\\
\frac{dy_{4}}{dX}&=&\frac{1}{\alpha(y_{1},y_{3},\widetilde{c}_{s})}\left(g+F(Y)\right),\label{adg10}
\end{eqnarray}
where
\begin{eqnarray}
\alpha(y_{1},y_{3},\widetilde{c}_{s})&=&\frac{3}{2}(y_{1}-y_{3})-\widetilde{c}_{s},\label{adg10b}\\
F(Y)&=&\left(2\widetilde{c}_{s}+\frac{1}{2}\right)y_{3}-\left(\widetilde{c}_{s}+\frac{1}{2}\right)y_{1}\nonumber\\
&&+\frac{1}{2}\left(\frac{3}{2}y_{1}^{2}+\frac{7}{2}y_{2}^{2}+y_{1}y_{3}-\frac{3}{2}y_{3}^{2}-6y_{2}y_{4}+3y_{4}^{2}\right).\label{adg10c}
\end{eqnarray}
The system (\ref{adg10})-(\ref{adg10c}) has a singularity as $\alpha=0$, meaning that
\begin{eqnarray*}
\frac{3}{2}\mathscr{J}u-\widetilde{c}_{s}=0\Rightarrow h=\frac{2}{3}\widetilde{c}_{s}.
\end{eqnarray*}
Out of this constant solution, the transformation $dX=\alpha dZ$ leads to the nonsingular system, \cite{ZhouT2010}
\begin{eqnarray}
\frac{dy_{1}}{dZ}&=&\alpha y_{2},\nonumber\\
\frac{dy_{2}}{dZ}&=&\alpha y_{3},\nonumber\\
\frac{dy_{3}}{dZ}&=&\alpha y_{4},\nonumber\\
\frac{dy_{4}}{dZ}&=&g+F(Y),\label{adg11}
\end{eqnarray}
We first discuss the equilibria of (\ref{adg11}). There are two cases:
\begin{itemize}
\item[(i)] $\alpha=0$. Then
$$y_{1}=y_{3}+\frac{2}{3}\left(c_{s}+\frac{1}{2}\right).$$ Furthermore, after some tedious algebra, the condition $g+F(Y)=0$ can be written as
\begin{eqnarray}
\frac{7}{2}x_{1}^{2}+\frac{3}{7}x_{2}^{2}+x_{3}^{2}+2 (g-g_{2}(c_{s}))=0,\label{adg12}
\end{eqnarray}
where
\begin{eqnarray}
x_{1}=y_{2}-\frac{6}{7}y_{4},\quad x_{2}=y_{4},\quad x_{3}=y_{3}+\frac{7}{3}\left(c_{s}+\frac{1}{2}\right),\label{adg13}
\end{eqnarray}
and
\begin{eqnarray}
g_{2}(c_{s}):=\frac{1}{3}\left(c_{s}+\frac{1}{2}\right)+\frac{52}{9}\left(c_{s}+\frac{1}{2}\right).\label{adg14}
\end{eqnarray}
Note that:
\begin{itemize}
\item If $g>g_{2}(c_{s})$, there are no points $(x_{1},x_{2},x_{3})$ satisfying (\ref{adg12}).
\item If $g=g_{2}(c_{s})$, the only point satisfying (\ref{adg12}) is $x_{1}=x_{2}=x_{3}=0$. Since the inverse of (\ref{adg13}) is
\begin{eqnarray}
y_{2}=x_{1}+\frac{6}{7}x_{2},\quad y_{4}=x_{2},\quad y_{3}=x_{3}-\frac{7}{3}\left(c_{s}+\frac{1}{2}\right),\label{adg15a}
\end{eqnarray}
and 
\begin{eqnarray}
y_{1}=y_{3}+\frac{2}{3}\left(c_{s}+\frac{1}{2}\right)=x_{3}-\frac{5}{3}\left(c_{s}+\frac{1}{2}\right),\label{adg15b}
\end{eqnarray} 
then the origin is transformed into
$$y_{1}=-\frac{5}{3}\left(c_{s}+\frac{1}{2}\right), y_{2}=0, y_{3}=-\frac{7}{3}\left(c_{s}+\frac{1}{2}\right), y_{4}=0.$$
\item If $g<g_{2}(c_{s})$, then we have the equilibria (\ref{adg15a}), (\ref{adg15b}), where $(x_{1},x_{2},x_{3})$ are the points of the ellipsoid (\ref{adg12}).
\end{itemize}
\item[(ii)] $\alpha\neq 0$. Then $y_{2}=y_{3}=y_{4}=0$ and the condition $g+F(Y)=0$ reads
\begin{eqnarray*}
y_{1}^{2}-\frac{4}{3}\left(c_{s}+1\right)y_{1}+\frac{4}{3}g=0,
\end{eqnarray*}
leading to the solutions
\begin{eqnarray}
y_{1}=y_{\pm}=\frac{2}{3}\left(\left(c_{s}+1\right)\pm \sqrt{\left(c_{s}+1\right)^{2}-{3}g}\right).\label{adg16}
\end{eqnarray}
Then, the solutions (\ref{adg16}) are real when 
\begin{eqnarray}
g\leq g_{1}(c_{s}):=\frac{}{3}\left(c_{s}+1\right).\label{adg17}
\end{eqnarray}
\end{itemize}
(Note that $g_{1}(c_{s})<g_{2}(c_{s})$ for all $c_{s}>0$.) We may summarize the previous analysis as follows. 
\begin{satz}
Let $c_{s}>0$ and $g_{1}(c_{s}), g_{2}(c_{s})$ be given respectively by (\ref{adg17}) and (\ref{adg14}). The following holds:
\begin{enumerate}
\item If $g<g_{1}(c_{s})$ then we have the equilibria 
\begin{eqnarray}
\widetilde{Y}&=&(y_{\pm},0,0,0),\label{adg17a}\\
Y^{*}&=&\left(x_{3}-\frac{5}{3}\left(c_{s}+\frac{1}{2}\right), x_{1}+\frac{6}{7}x_{2},x_{3}-\frac{7}{3}\left(c_{s}+\frac{1}{2}\right),x_{2}\right),\label{adg17b} 
\end{eqnarray}
where $y_{\pm}$ is given by (\ref{adg16}) and $(x_{1},x_{2},x_{3})$ are the points of the ellipsoid (\ref{adg12}).
\item If $g=g_{1}(c_{s})$, then (\ref{adg17b}) holds and (\ref{adg17a}) becomes
$$\widetilde{Y}=(y_{1},0,0,0),\quad y_{1}=y_{\pm}=\frac{2}{3}\left(c_{s}+1\right).$$
\item If $g_{1}(c_{s})<g<g_{2}(c_{s})$, then we just have the equilibria (\ref{adg17b}).
\item If $g=g_{2}(c_{s})$, then the equilibria (\ref{adg17b}) reduces to
$$Y^{*}=\left(-\frac{5}{3}\left(c_{s}+\frac{1}{2}\right), 0,-\frac{7}{3}\left(c_{s}+\frac{1}{2}\right),0\right),$$ corresponding to $x_{1}=x_{2}=x_{3}=0$.
\item If $g>g_{2}(c_{s})$, then there are no equilibria.
\end{enumerate}
\end{satz}
\subsection{The case $g\leq g_{1}(c_{s})$} 
For $c_{s}>0$ fixed, the present paper is focused on the existence of homolcinic and periodic orbits around the equilibria (\ref{adg17a}). The first group will correspond to traveling wave (TW) solutions of (\ref{adg8}) or (\ref{adg8b}) approaching  $y_{+}$ (or $y_{-}$) when $t\rightarrow\pm\infty$. The second group is associated to the existence of periodic traveling wave (PTW) solutions. Note in particular that when $g=0$ then
$$y_{+}=\frac{}{3}\left(c_{s}+1\right),\quad y_{-}=0,$$ and those orbits homo-clinic to $y_{-}$ at infinity are identified as solitary wave solutions, with the meaning of localized traveling wave solutions.

\begin{figure}[htbp]
\centering
{\includegraphics[width=1\textwidth]{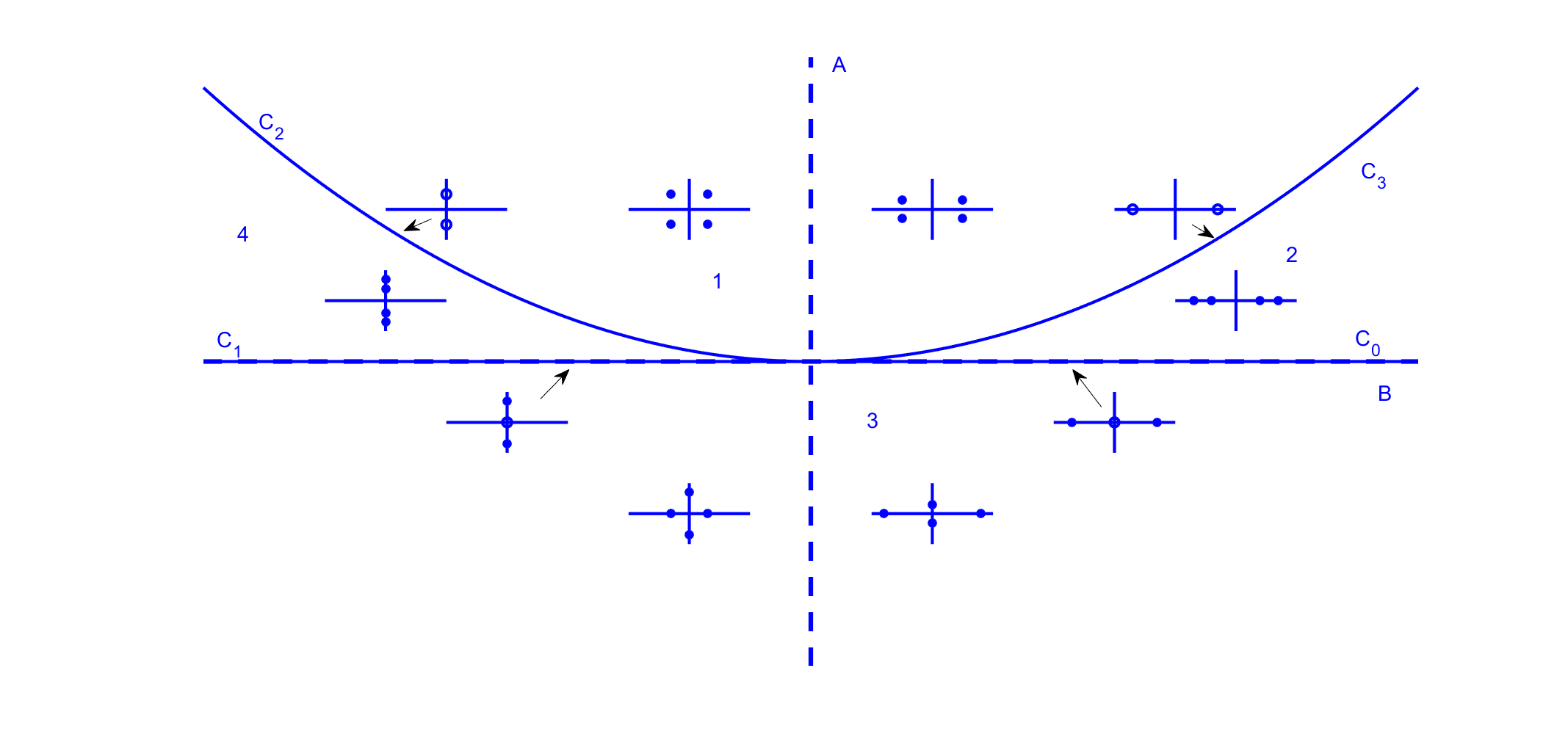}}
\caption{Linearization at the origin of (\ref{adg19}) (cf. \cite[Figure 1]{Champ}): Regions $1$ to $4$ in the $(b,a)$-plane, delimited by the bifurcation curves $\mathbb{C}_{0}$ to $\mathbb{C}_{3}$ given by (\ref{bifurcurv}), and schematic representation of the position in the complex plane of the eigenvalues of $L$ for each curve and region. (Dot: simple root, larger dot: double root.)}
\label{ADG0}
\end{figure}

To this end, our approach will first look for a reduction of (\ref{adg9c}) to the case $g=0$. Let $y_{\pm}$ denote any of the equilibria (\ref{adg16}). Note that if $u$ is a solution of  (\ref{adg9c}), then $\widetilde{u}=u-y_{\pm}$ satisfies
\begin{eqnarray}
\left(\frac{3}{2}\mathscr{J}\widetilde{u}+\frac{3y_{\pm}-1}{2}-{c}_{s}\right)\widetilde{u}''''+\left(2{c}_{s}+\frac{3-5y_{\pm}}{2}\right)\widetilde{u}''+\left(-{c}_{s}+\frac{3y_{\pm}-2}{2}\right)\widetilde{u}&&\nonumber\\
=f(\widetilde{u},\widetilde{u}',\widetilde{u}'',\widetilde{u}''')&&\nonumber\\
=\frac{1}{2}\left(-\frac{3}{2}\widetilde{u}^{2}+\frac{5}{2}(\widetilde{u}')^{2}+5\widetilde{u}\widetilde{u}''-\frac{9}{2}(\widetilde{u}'')^{2}-6\widetilde{u}'\widetilde{u}'''+3(\widetilde{u}''')^{2}\right),&&\label{adg18}
\end{eqnarray}
As before, (\ref{adg18}) can be written as a first-order system for
$U=(U_{1},U_{2},U_{3},U_{4})^{T}=(\widetilde{u},\widetilde{u}',\widetilde{u}'',\widetilde{u}''')^{T}$ as
\begin{eqnarray}
\frac{dU}{dZ}=L(c_{s},y_{\pm})U+G(U),\label{adg19}
\end{eqnarray}
where
\begin{eqnarray*}
L&=&L(c_{s},y_{\pm})=\begin{pmatrix}0&\frac{3y_{\pm}-1}{2}-{c}_{s}&0&0\\
0&0&\frac{3y_{\pm}-1}{2}-{c}_{s}&0\\
0&0&0&\frac{3y_{\pm}-1}{2}-{c}_{s}\\
{c}_{s}+\frac{2-3y_{\pm}}{2}&0&-(2{c}_{s}+\frac{3-5y_{\pm}}{2})
\end{pmatrix},\\
G(U)&=&\begin{pmatrix}\frac{3}{2}U_{2}(U_{1}-U_{3})\\\frac{3}{2}U_{3}(U_{1}-U_{3})\\\frac{3}{2}U_{4}(U_{1}-U_{3})\\\widetilde{F}(U)\end{pmatrix},\\
\widetilde{F}(U)&=&-\left(\frac{9}{4}U_{3}^{2}-\frac{5}{2}U_{1}U_{3}\right)-\frac{3}{4}U_{1}^{2}+\frac{5}{4}U_{2}^{2}-{3}U_{2}U_{4}-\frac{3}{2}U_{4}^{2}.
\end{eqnarray*}
The structure of solutions of (\ref{adg8}) or (\ref{adg8b}) around $y_{\pm}$ can be analyzed from that of solutions of (\ref{adg19}) around $U=0$. Since $G(0)=G'(0)=0$, then, if we first linearize, the characteristic polynomial of $L$ is of the form $z^{4}-bz^{2}+a$, where
\begin{eqnarray}
&&a=a_{\pm}=-\beta_{\pm}\alpha_{\pm}^{3},\quad b=b_{\pm}=-\alpha_{\pm}\gamma_{\pm},\nonumber\\
&&\alpha_{\pm}=\frac{3y_{\pm}-1}{2}-{c}_{s},\quad \beta_{\pm}={c}_{s}+\frac{2-3y_{\pm}}{2},\quad \gamma_{\pm}=2{c}_{s}+\frac{3-5y_{\pm}}{2},\label{adg20}
\end{eqnarray}
We note that (\ref{adg19}) is reversible with respect to the transformation
$$S:(U_{1},U_{2},U_{3},U_{4})\mapsto (U_{1},-U_{2},U_{3},-U_{4}),$$ in the sense that, \cite{HaragusI}
$$SLU=-LSU,\quad SG(U)=-G(SU).$$
This implies that the linear dynamics around $U=0$ follows Figure \ref{ADG0}, \cite{Champ}, in the $(b,a)$ plane. Four regions of different dynamics are determined from the four curves
\begin{eqnarray}
\mathbb{C}_{0}&=&\{(b,a) / a=0, b>0\},\nonumber\\
\mathbb{C}_{1}&=&\{(b,a) / a=0, b<0\},\nonumber\\
\mathbb{C}_{2}&=&\{(b,a) / a>0, b=-2\sqrt{a}\},\nonumber\\
\mathbb{C}_{3}&=&\{(b,a) / a>0, b=2\sqrt{a}\}.\label{bifurcurv}
\end{eqnarray}

On the curve $\mathbb{C}_{0}$ there are two zero eigenvalues and two real; on $\mathbb{C}_{1}$ there are two zero eigenvalues and two imaginary; on $\mathbb{C}_{2}$ there is a double complex conjugate pair of imaginary eigenvalues $\pm i\sqrt{|b|/2}$, and on $\mathbb{C}_{3}$ there are two double real eigenvalues $\pm\sqrt{b/2}$, symmetric with respect to the imaginary axis. The spectrum of $L$ in the regions delimited by the four curves is represented schematically in Figure \ref{ADG0}.

The existence of different types of orbits in each region and with respect to each equilibrium in (\ref{adg16}) is analyzed using two approaches. The first one is theoretical and is based on the Normal Form theory and the theory of reversible bifurcations, cf. e.~g. \cite{HaragusI}, as well as the study of the corresponding systems derived from a center manifold reduction. But the discussion will be also strongly supported by numerical simulations. To this end, a numerical procedure to approximate solutions of (\ref{adg18}) is implemented. This is described in Appendix \ref{appA}.
\subsection{Some previous lemmas}
\begin{figure}[htbp]
\centering
\centering
\subfigure[]
{\includegraphics[width=0.41\columnwidth]{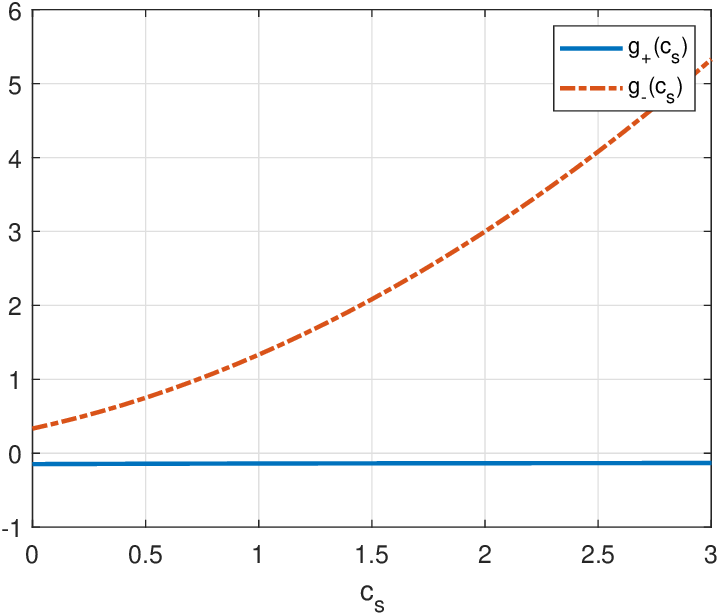}}
\subfigure[]
{\includegraphics[width=0.45\columnwidth]{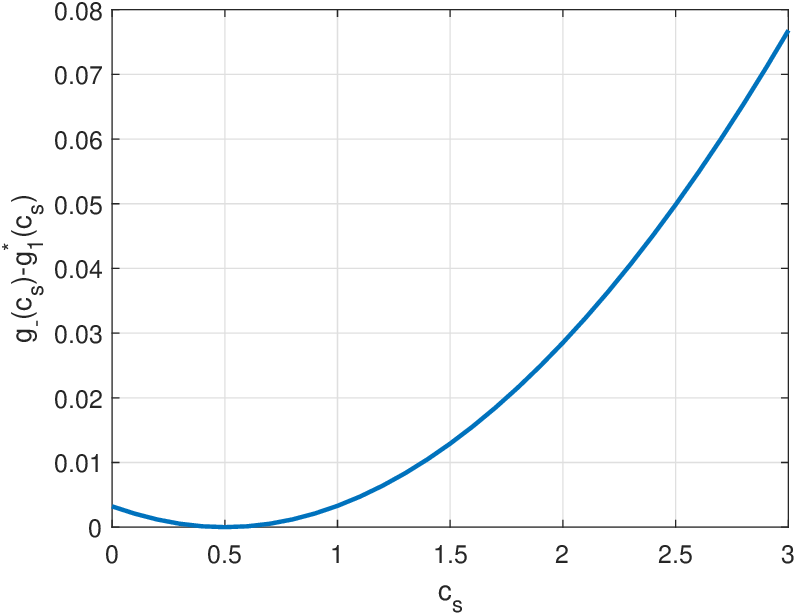}}
\subfigure[]
{\includegraphics[width=0.45\columnwidth]{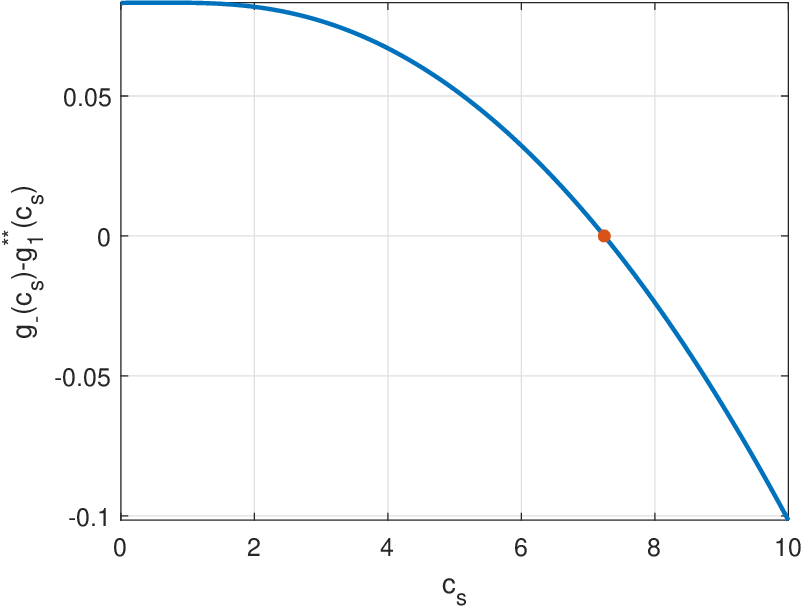}}
\subfigure[]
{\includegraphics[width=0.45\columnwidth]{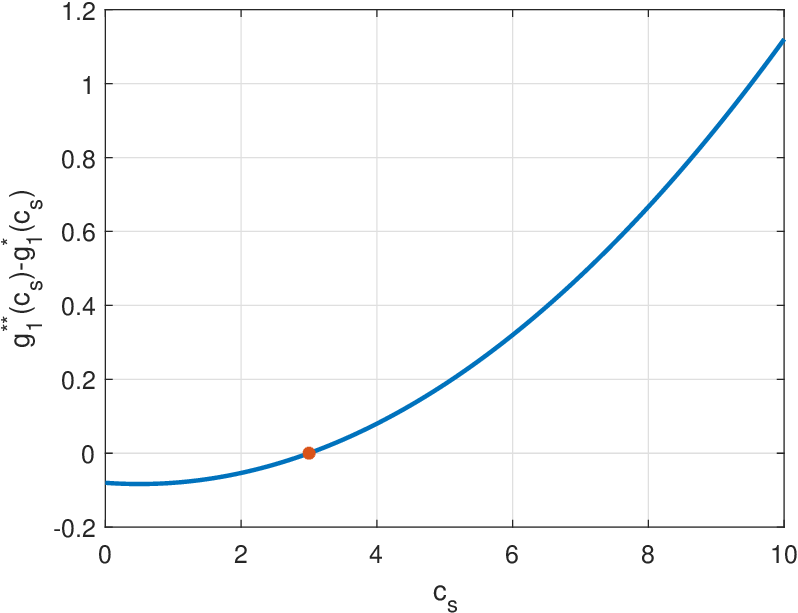}}
\caption{Illustration of Lemma \ref{adg_lemma2}. (a), (b) Case (1); (c) Case (2), where the remarked dot corresponds to $c_{s}= 3(1+\sqrt{2})$; (d) Case (c), where the remarked dot corresponds to $c_{s}= 3$}
\label{ADG00}
\end{figure}

The analysis below will require some previous results concerning the coefficients (\ref{adg20}) and new relevant values of the constant $g$.
We note first that when $g\leq g_{1}(c_{s})$ then $0\leq y_{-}\leq y_{+}$. Furthermore, it is not hard to prove the following result.
\begin{lemma}
\label{adg_lemma1}
Let $c_{s}>0$ and assume that $g\leq g_{1}(c_{s})$. we define
\begin{eqnarray}
g_{1}^{*}(c_{s}):=g_{1}(c_{s})-\frac{1}{75}\left(c_{s}-\frac{1}{2}\right)^{2},\quad
g_{1}^{**}(c_{s}):=g_{1}(c_{s})-\frac{1}{12}.\label{newg}
\end{eqnarray}
Then it holds that:
\begin{itemize}
\item[(1)] $\alpha_{+}>0, \beta_{+}\leq 0, \beta_{-}\geq 0$.
\item[(2)] $\alpha_{-}\geq 0\Leftrightarrow g\geq g_{1}^{**}(c_{s})$.
\item[(3)] $\gamma_{+}>0\Leftrightarrow c_{s}>\frac{1}{2}$ and $g>g_{1}^{*}(c_{s})$.
\item[(4)] $\gamma_{-}>0\Leftrightarrow c_{s}\geq\frac{1}{2}$ and $g<g_{1}^{*}(c_{s})$ when $c_{s}<\frac{1}{2}$.
\end{itemize}
\end{lemma}
In addition, for each case $y_{+}, y_{-}$, we will 
need to identify the position of the roots of characteristic polynomial of $L$ in the $(b,a)$ plane from the sign of the coefficients $a_{\pm}, b_{\pm}$ and $\Delta_{\pm}=b_{\pm}^{2}-4a_{\pm}$. After some computations, we observe that
\begin{eqnarray*}
\Delta_{\pm}&=&\alpha_{\pm}^{2}S_{\pm},\\
S_{\pm}&=&\gamma_{\pm}^{2}+4\alpha_{\pm}\beta_{\pm}=-\frac{11}{9}(\beta_{\pm}-\delta_{+})(\beta_{\pm}-\delta_{-}),\\
\delta_{\pm}&=&\frac{1}{22}\left(\left(10c_{s}+{13}\right)\pm\sqrt{\left(10c_{s}+{13}\right)^{2}+44\left(c_{s}-\frac{1}{2}\right)^{2}}\right),
\end{eqnarray*}
with $\delta_{-}<0<\delta_{+}$ for $c_{s}>0$. This leads to the following properties:
\begin{lemma}
\label{adg_lemma1b}
Let $c_{s}>0$ and assume that $g\leq g_{1}(c_{s})$. We define
\begin{eqnarray}
g_{-}(c_{s}):=g_{1}(c_{s})-\frac{1}{3}\delta_{-}^{2},\quad
g_{+}(c_{s}):=g_{1}(c_{s})-\frac{1}{3}\delta_{+}^{2}.\label{newg2}
\end{eqnarray}
Then:
\begin{itemize}
\item[(i)] $\Delta_{+}>0\Leftrightarrow g>g_{-}(c_{s})$.
\item[(ii)] If $g_{1}^{**}(c_{s})<g<g_{1}(c_{s})$ then $\Delta_{-}>0$. When $g<g_{1}^{**}(c_{s})$ then 
$$\Delta_{-}>0\Leftrightarrow g>g_{+}(c_{s}).$$
\end{itemize}
\end{lemma}
A final result compares the limiting values of $g$ appearing above and in Lemma \ref{adg_lemma1} and it is illustrated in Figure \ref{ADG00}.
\begin{lemma}
\label{adg_lemma2}
Let $c_{s}>0$ and assume that $g\leq g_{1}(c_{s})$.
\begin{itemize}
\item[(1)] $g_{+}(c_{s})<0<g_{-}(c_{s})$ and $g_{1}^{*}(c_{s}), g_{1}^{**}(c_{s})>0$.
\item[(2)] $g_{1}^{*}(c_{s})<g_{-}(c_{s})$ and $g_{1}^{**}(c_{s})\leq g_{-}(c_{s})\Leftrightarrow c_{s}\leq 3(1+\sqrt{2})$.
\item[(3)] $g_{1}^{*}(c_{s})\leq g_{1}^{**}(c_{s})\Leftrightarrow c_{s}\geq 3$,
\end{itemize}
where $g_{1}^{*}(c_{s}), g_{1}^{**}(c_{s})$ and $g_{+}(c_{s}), g_{-}(c_{s})$ are given by (\ref{newg}) and (\ref{newg2}), respectively.
\end{lemma}

\section{Existence of non-periodic traveling wave solutions}
\label{sec3}
Making use of Lemmas \ref{adg_lemma1}-\ref{adg_lemma2}, we will study the dynamics of 
(\ref{adg19}) with respect to each of the equilibria (\ref{adg16}) 
\subsection{Dynamics around the equilibrium $y_{+}$}
The analysis for the equilibrium $y_{+}$ leads to the following situations:
\begin{itemize}
\item[(i)] Assume that $c_{s}>\frac{1}{2}$. In that case, $a_{+}>0$ and we have three possibilities:
\begin{itemize}
\item[(i1)] If $g<g_{1}^{*}(c_{s})$, then $b_{+}>0$ and $\Delta_{+}<0$. The linearization matrix $L(c_{s},y_{+})$ has two symmetric couples of complex conjugate eigenvalues (region 1, right of Figure \ref{ADG0}). Here the dynamics can be analyzed following the approach made in, e.~g. \cite{ChampT}. The numerical computations suggest the generation of classical traveling waves of depression with a nonmonotone behaviour. For $g$ close to $g_{1}^{*}(c_{s})$, the profiles are smooth (cf. Figure \ref{ADG4a}), developing some kind of peak at the point of maximum negative excursion as $g$ separates from $g_{1}^{*}(c_{s})$.
\item[(i2)] When $g_{1}^{*}(c_{s})\leq g<g_{-}(c_{s})$, then $b_{+}<0$ and $\Delta_{+}<0$. The spectrum of the linearization corresponds to region 1, left, of Figure \ref{ADG0}. In Figure \ref{ADG5a}, it is shown the persistence of the generation of classical traveling waves (of depression) with nonmonotone asymptotic behaviour. Note that as $g\rightarrow g_{-}(c_{s})$ then $\Delta_{+}\rightarrow 0$ and we fall into the curve $\mathbb{C}_{2}$, where the spectrum of the linearization changes to a doble complex conjugate pair of imaginary eigenvalues. This seems to be reflected by a more oscillatory approach towards the equlibrium $y_{+}$ as $Z\rightarrow\pm\infty$, \cite{DougalisDM2012,Champ}, see Figure \ref{ADG5b}.
\item[(i3)] When $g_{-}(c_{s})\leq g<g_{1}(c_{s})$, then $b_{+}<0$ and $\Delta_{+}>0$ (region 4 of Figure \ref{ADG0}, where the spectrum has four imaginary eigenvalues, two by two conjugate). The behaviour of the orbits as $g\rightarrow g_{-}$, described above, and the numerical computations (cf. Figure \ref{ADG6a}) suggest that crossing from region 1 (left) to region 4 forms a bifurcation that may change the dynamics from the generation of nonmonotone traveling waves to periodic orbits in the form of periodic traveling waves,  \cite{IK}. As $g$ approaches $g_{1}(c_{s})$, we are approximating to the curve $\mathbb{C}_{1}$, where the spectrum consists of two zero eigenvalues and two imaginary. The computations do not detect any modification in the nature of the generated profiles, although the amplitude of the emerging periodic traveling waves seems to decrease to zero and the solution tends to the constant state $y_{+}$.

\begin{figure}[htbp]
\centering
\centering
\subfigure[]
{\includegraphics[width=0.45\columnwidth]{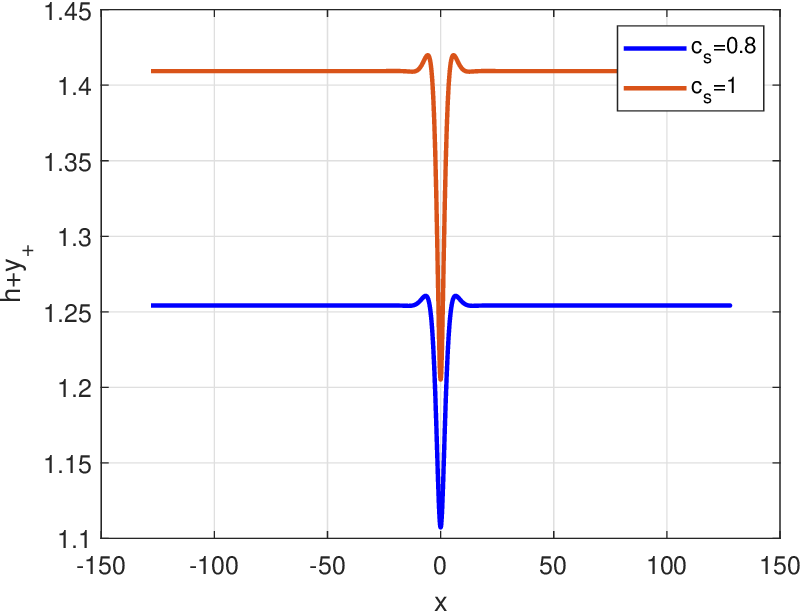}}
\subfigure[]
{\includegraphics[width=0.45\columnwidth]{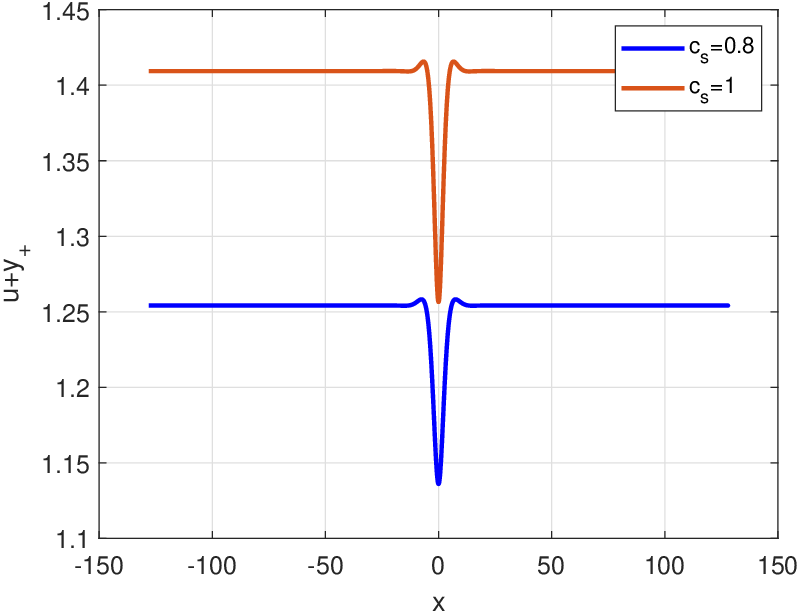}}
\subfigure[]
{\includegraphics[width=0.45\columnwidth]{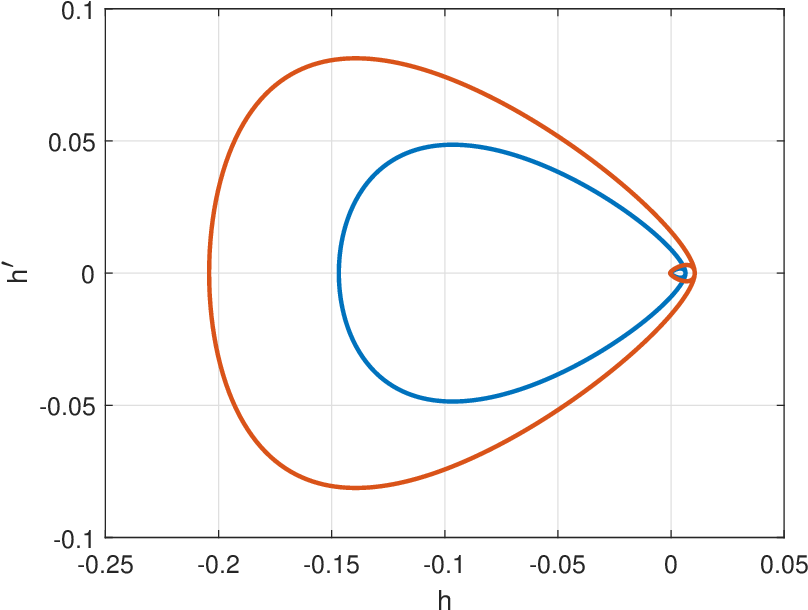}}
\subfigure[]
{\includegraphics[width=0.45\columnwidth]{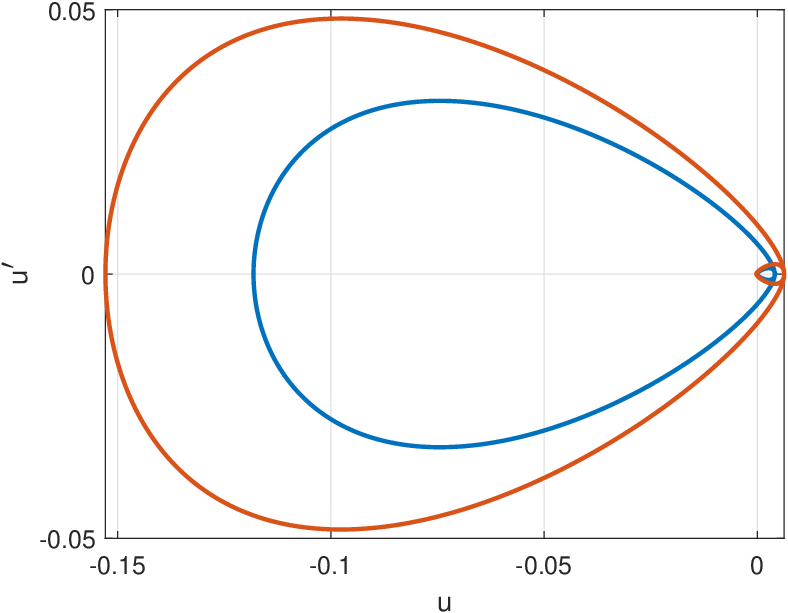}}
\caption{Approximations to traveling-wave solutions of (\ref{adg9}) for two values of $c_{s}$. Case $g=g_{1}^{*}(c_{s})-0.001$. (a) $\widetilde{h}+y_{+}$ profile; (b) $\widetilde{u}+y_{+}$ profile; (c) $\widetilde{h}$ phase plot; (d) $\widetilde{u}$ phase plot.}
\label{ADG4a}
\end{figure}
\begin{figure}[htbp]
\centering
\centering
\subfigure[]
{\includegraphics[width=0.45\columnwidth]{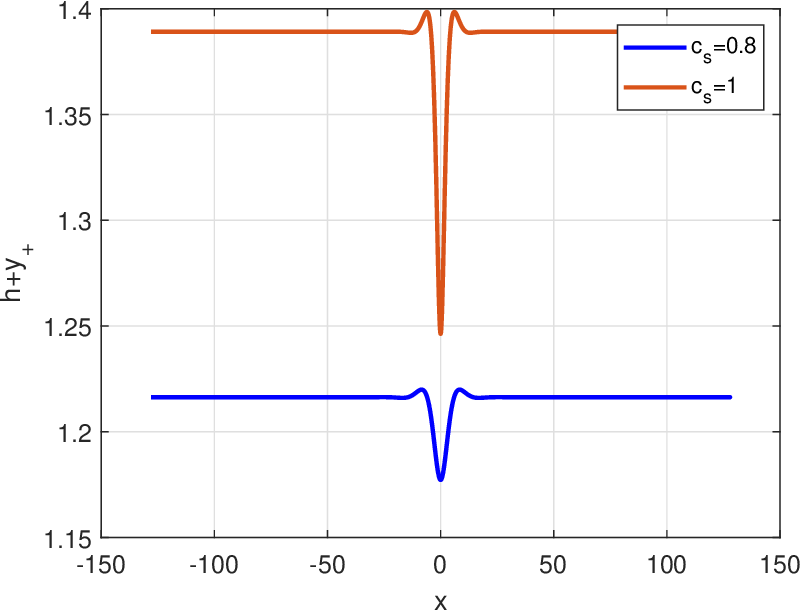}}
\subfigure[]
{\includegraphics[width=0.45\columnwidth]{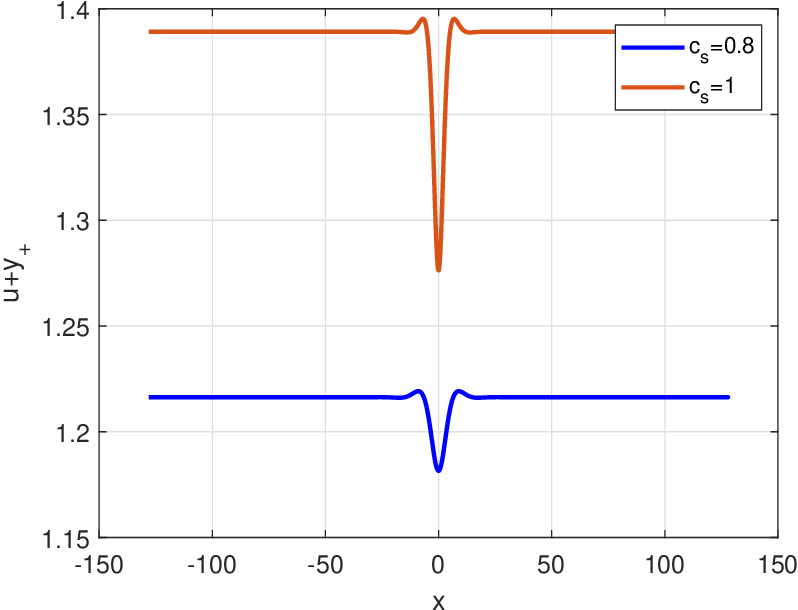}}
\subfigure[]
{\includegraphics[width=0.45\columnwidth]{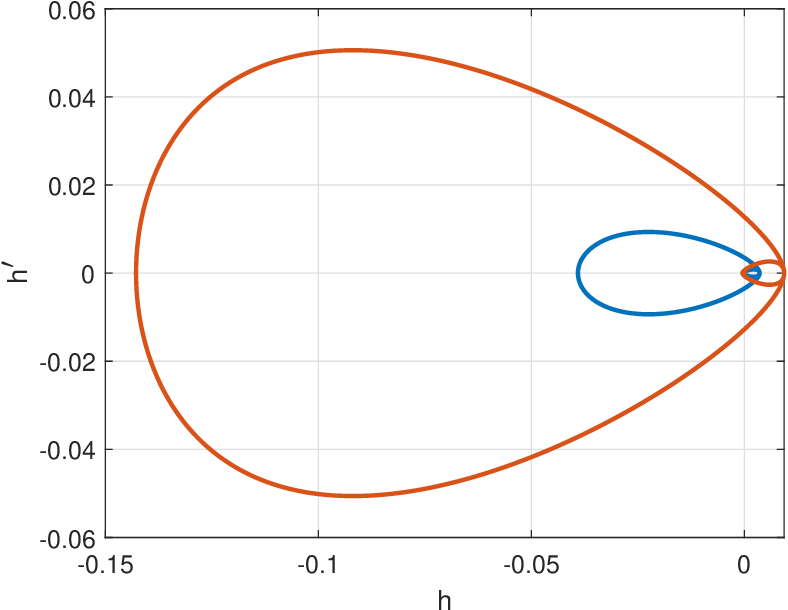}}
\subfigure[]
{\includegraphics[width=0.45\columnwidth]{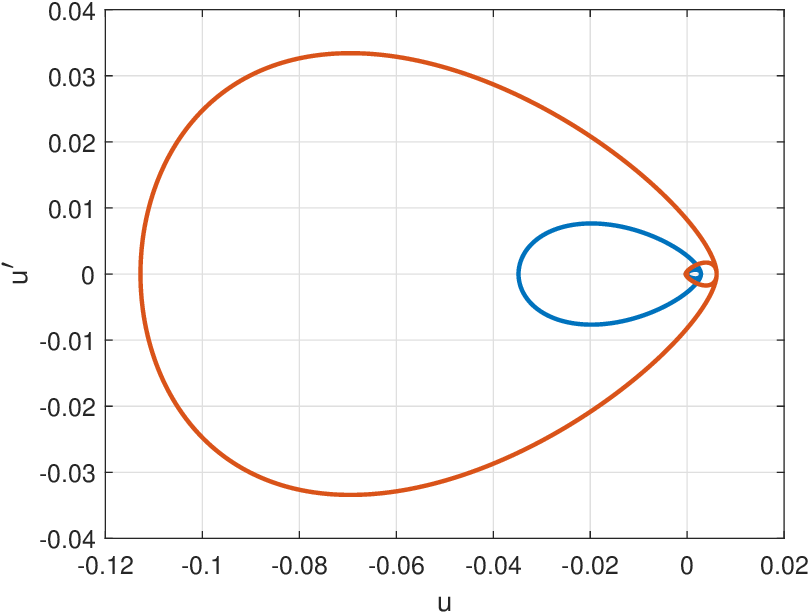}}
\caption{Approximations to traveling-wave solutions of (\ref{adg9}) for two values of $c_{s}$. Case $g=g_{1}^{*}(c_{s})+0.001$. (a) $\widetilde{h}+y_{+}$ profile; (b) $\widetilde{u}+y_{+}$ profile; (c) $\widetilde{h}$ phase plot; (d) $\widetilde{u}$ phase plot.}
\label{ADG5a}
\end{figure}
\begin{figure}[htbp]
\centering
\centering
\subfigure[]
{\includegraphics[width=0.45\columnwidth]{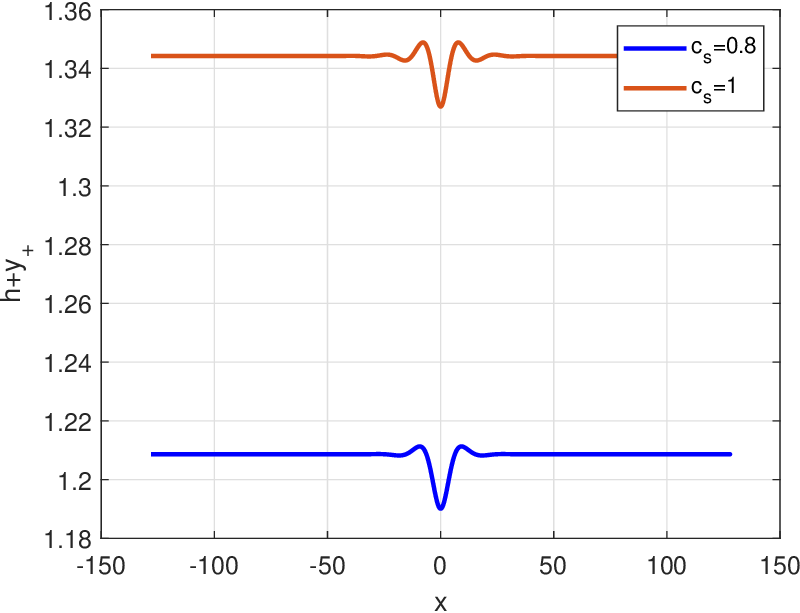}}
\subfigure[]
{\includegraphics[width=0.45\columnwidth]{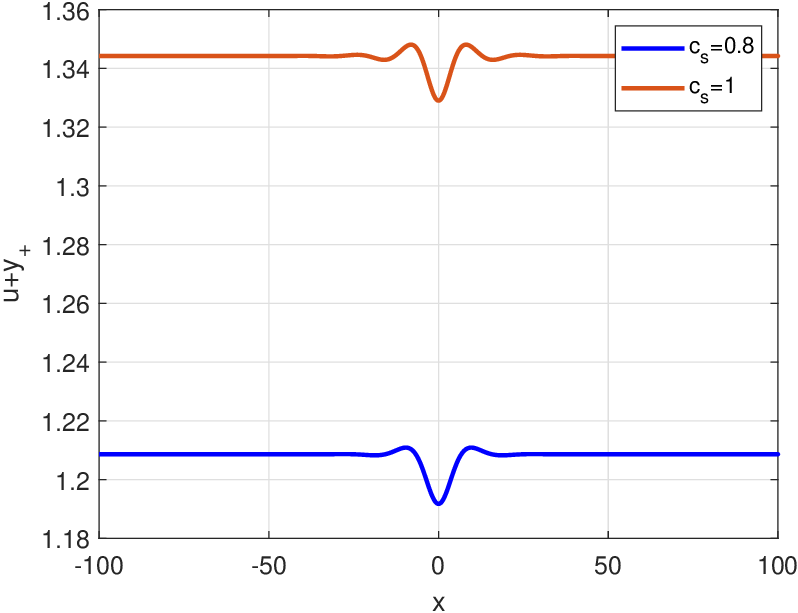}}
\subfigure[]
{\includegraphics[width=0.45\columnwidth]{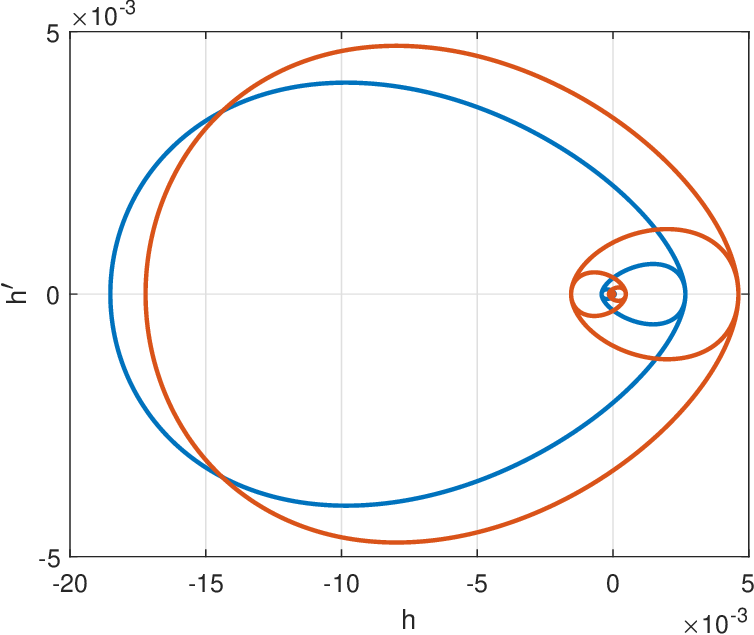}}
\subfigure[]
{\includegraphics[width=0.45\columnwidth]{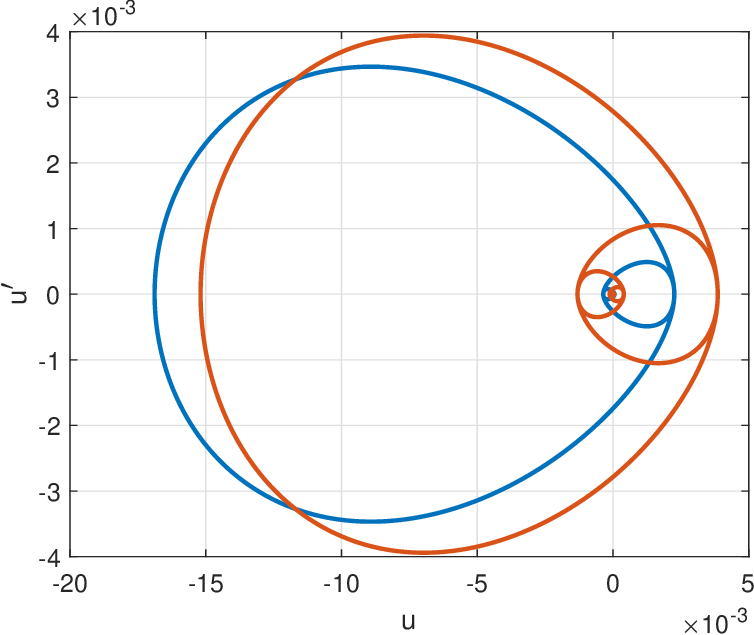}}
\caption{Approximations to traveling-wave solutions of (\ref{adg9}) for two values of $c_{s}$. Case $g=g_{-}(c_{s})-10^{-4}$. (a) $\widetilde{h}+y_{+}$ profile; (b) $\widetilde{u}+y_{+}$ profile; (c) $\widetilde{h}$ phase plot; (d) $\widetilde{u}$ phase plot.}
\label{ADG5b}
\end{figure}
\begin{figure}[htbp]
\centering
\centering
\subfigure[]
{\includegraphics[width=0.45\columnwidth]{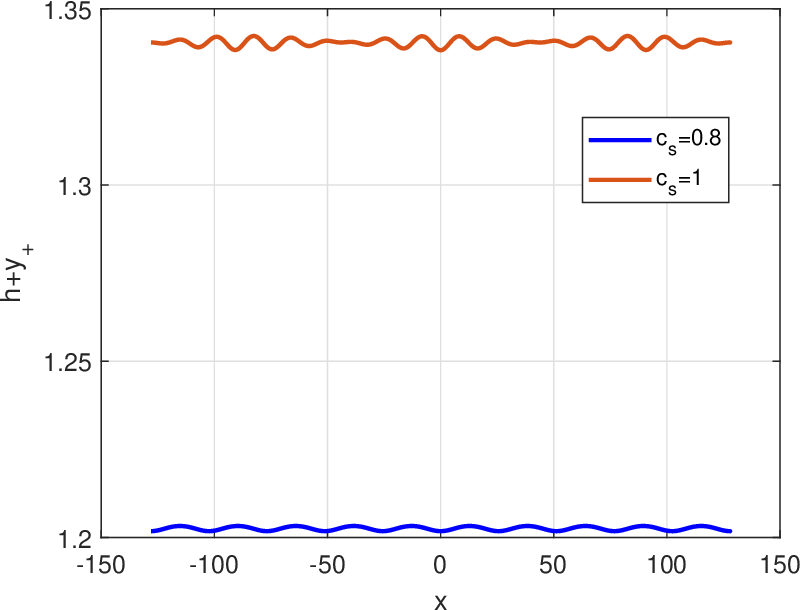}}
\subfigure[]
{\includegraphics[width=0.45\columnwidth]{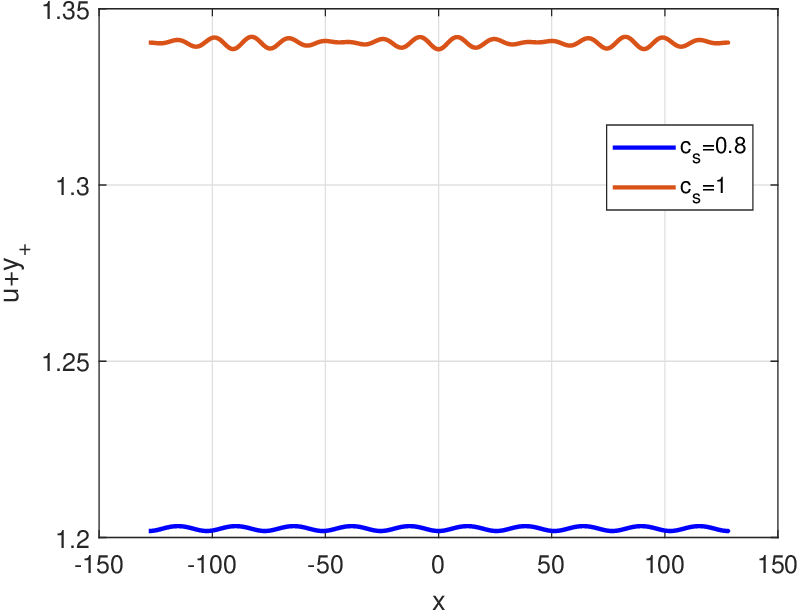}}
\subfigure[]
{\includegraphics[width=0.45\columnwidth]{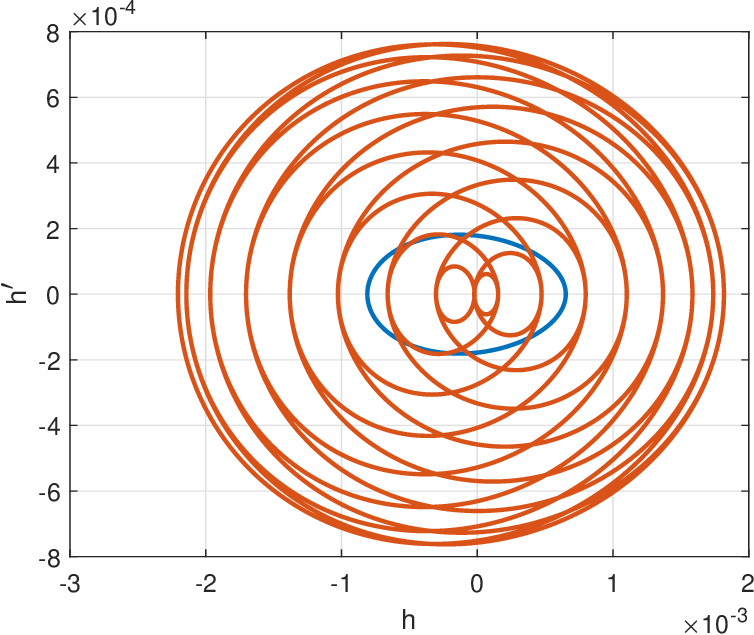}}
\subfigure[]
{\includegraphics[width=0.45\columnwidth]{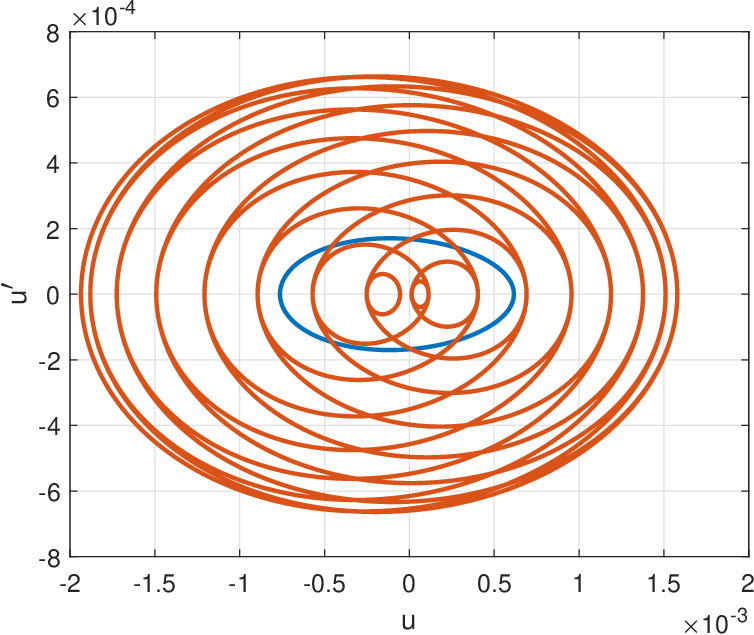}}
\caption{AApproximations to traveling-wave solutions of (\ref{adg9}) for two values of $c_{s}$. Case $g=g_{-}(c_{s})+10^{-6}$. (a) $\widetilde{h}+y_{+}$ profile; (b) $\widetilde{u}+y_{+}$ profile; (c) $\widetilde{h}$ phase plot; (d) $\widetilde{u}$ phase plot.}
\label{ADG6a}
\end{figure}
\end{itemize}
\item[(ii)] Assume that $c_{s}<\frac{1}{2}$. In that case, $a_{+}>0$, $b_{+}>0$ (since $\gamma_{+}<0$) and we may have two possibilities:
\begin{itemize}
\item[(ii1)] If $g<g_{-}(c_{s})$ then $\Delta_{+}<0$ and this corresponds to region 1, right of Figure \ref{ADG0}: the linearization matrix $L(c_{s},y_{+})$  has four (symmetric) complex eigenvalues. The results are then similar to those obtained in (i1), where the homo-clinic orbits correspond to nonmonotone classical traveling waves. In this case, the profiles seem to be more regular as $g$ approaches $g_{-}(c_{s})$ (cf. Figures \ref{ADG7a} and \ref{ADG7b}).
\item [(ii2)]  If $g_{-}(c_{s})\leq g<g_{1}(c_{s})$ then $\Delta_{+}>0$ and we fall into region 2 of Figure \ref{ADG0}. When crossing from region 1 to region 2 and $g$ goes from $g_{-}$ to $g_{1}$ (approaching the curve $\mathbb{C}_{0}$) there seems to be a change of type of homo-clinic orbit (cf Figures \ref{ADG8a} and \ref{ADG8b}) with the generation of classical traveling waves with a monotone asymptotic behaviour towards $y_{+}$, \cite{IK}.
\end{itemize}
\begin{figure}[htbp]
\centering
\centering
\subfigure[]
{\includegraphics[width=0.45\columnwidth]{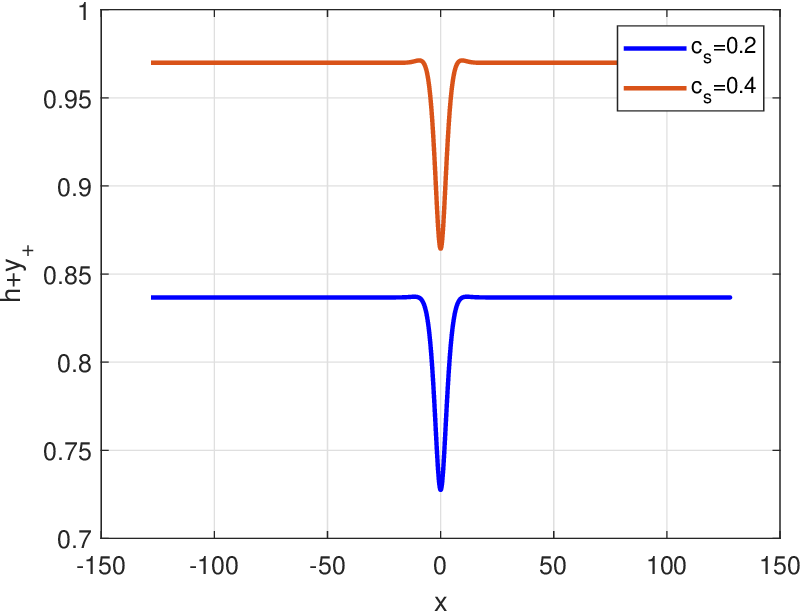}}
\subfigure[]
{\includegraphics[width=0.45\columnwidth]{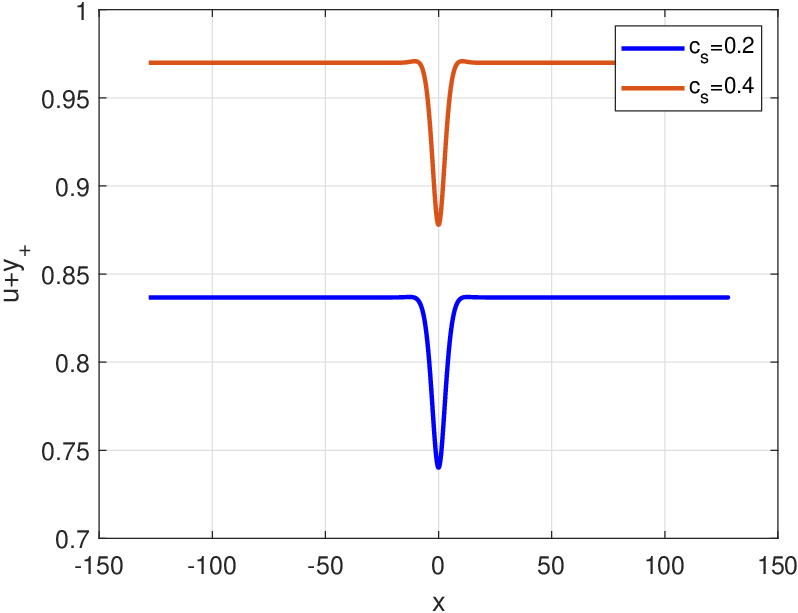}}
\subfigure[]
{\includegraphics[width=0.45\columnwidth]{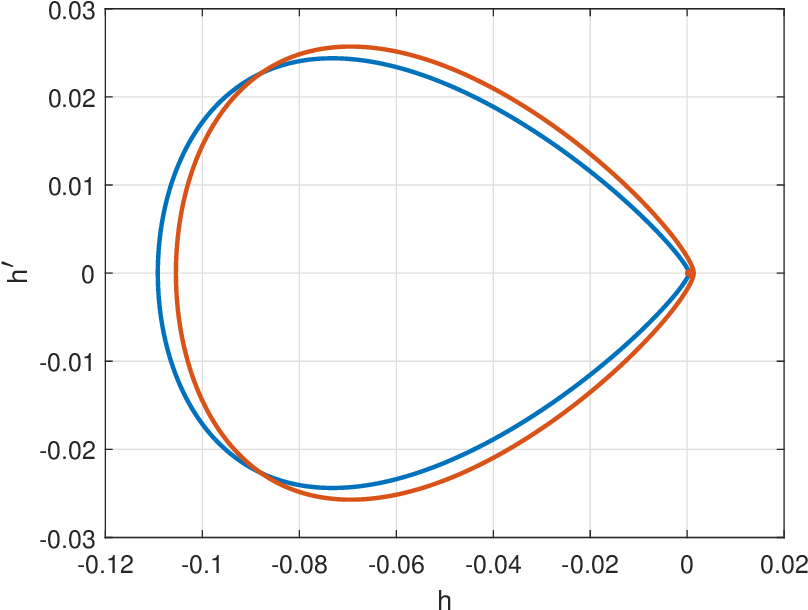}}
\subfigure[]
{\includegraphics[width=0.45\columnwidth]{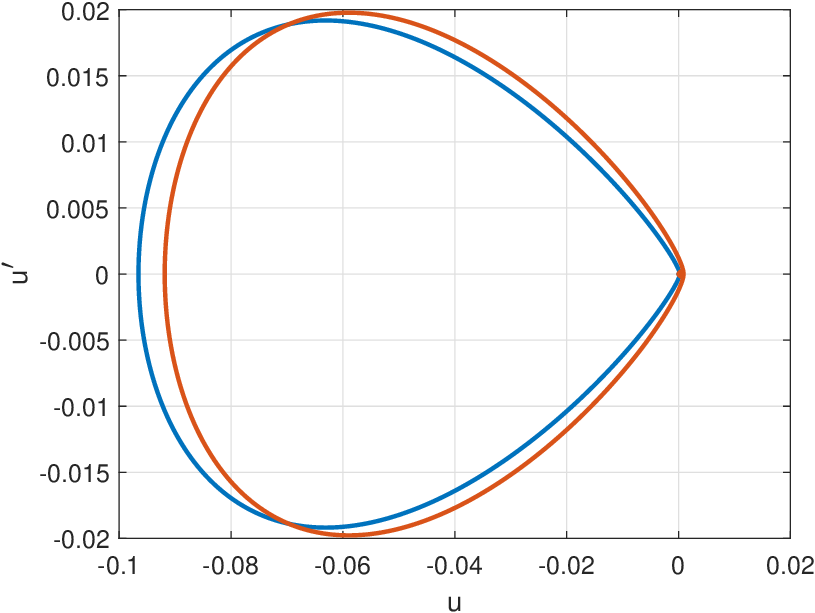}}
\caption{Approximations to traveling-wave solutions of (\ref{adg9}) for two values of $c_{s}$. Case $g=g_{-}(c_{s})-0.001$. (a) $\widetilde{h}+y_{+}$ profile; (b) $\widetilde{u}+y_{+}$ profile; (c) $\widetilde{h}$ phase plot; (d) $\widetilde{u}$ phase plot.}
\label{ADG7a}
\end{figure}
\begin{figure}[htbp]
\centering
\centering
\subfigure[]
{\includegraphics[width=0.45\columnwidth]{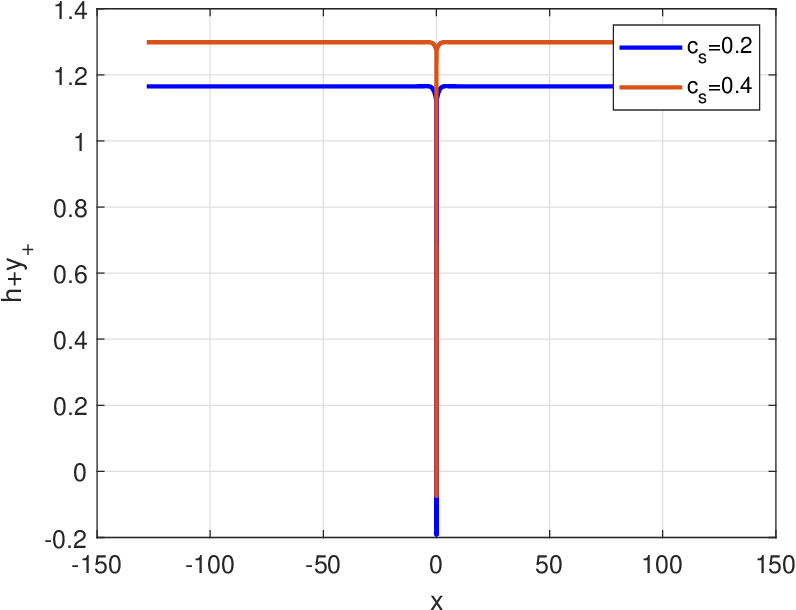}}
\subfigure[]
{\includegraphics[width=0.45\columnwidth]{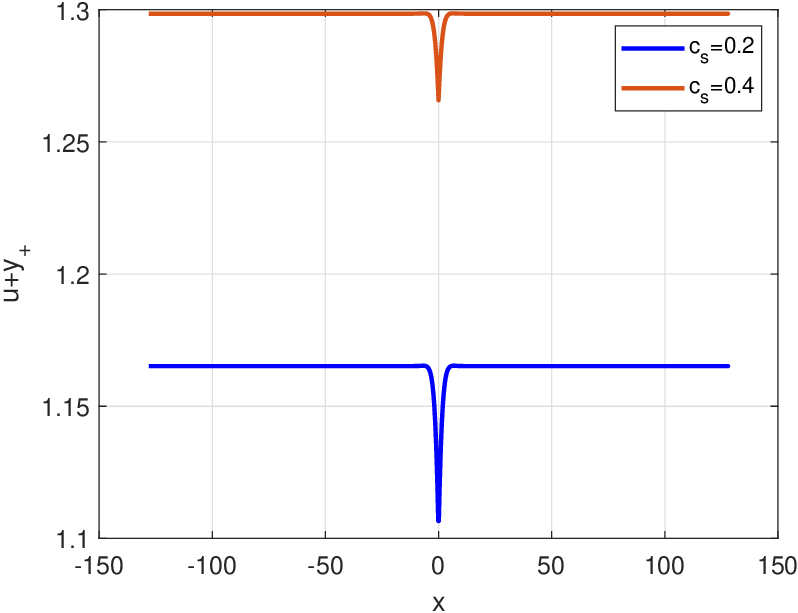}}
\subfigure[]
{\includegraphics[width=0.45\columnwidth]{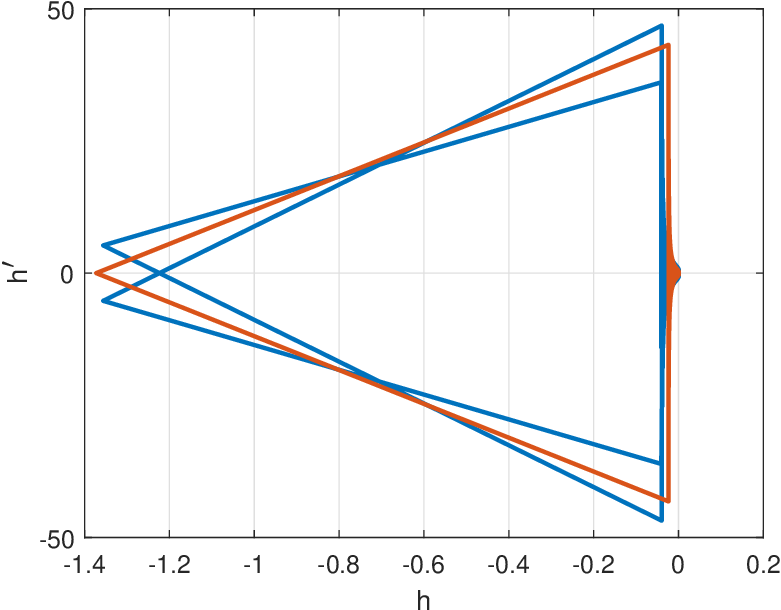}}
\subfigure[]
{\includegraphics[width=0.45\columnwidth]{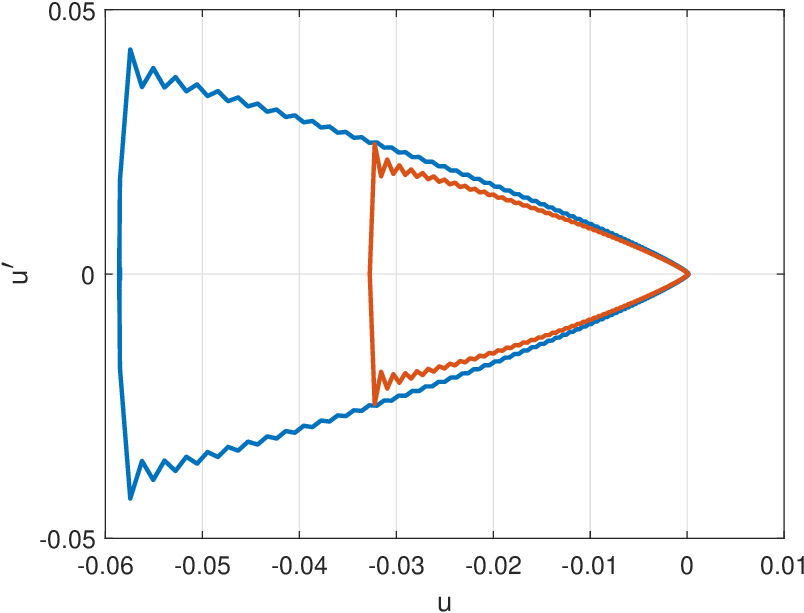}}
\caption{Approximations to traveling-wave solutions of (\ref{adg9}) for two values of $c_{s}$. Case $g=g_{-}(c_{s})-0.001$. (a) $\widetilde{h}+y_{+}$ profile; (b) $\widetilde{u}+y_{+}$ profile; (c) $\widetilde{h}$ phase plot; (d) $\widetilde{u}$ phase plot.}
\label{ADG7b}
\end{figure}
\begin{figure}[htbp]
\centering
\centering
\subfigure[]
{\includegraphics[width=0.45\columnwidth]{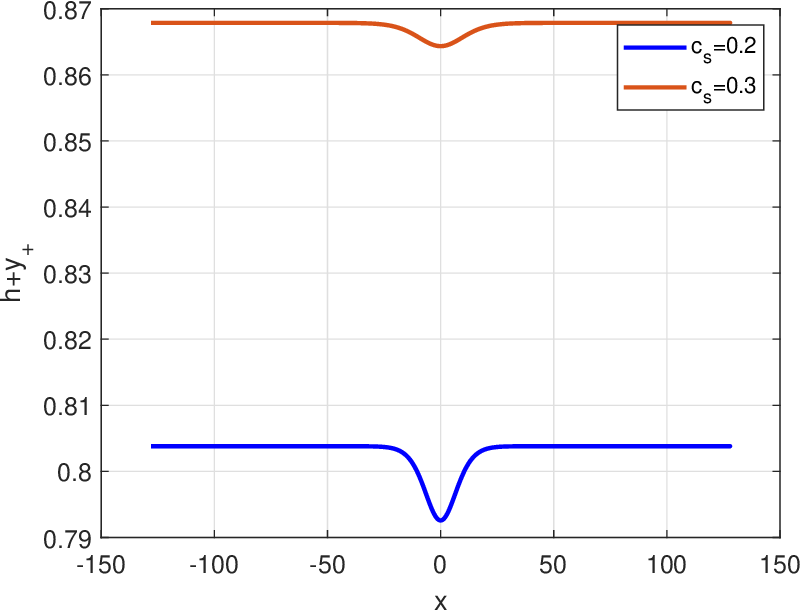}}
\subfigure[]
{\includegraphics[width=0.45\columnwidth]{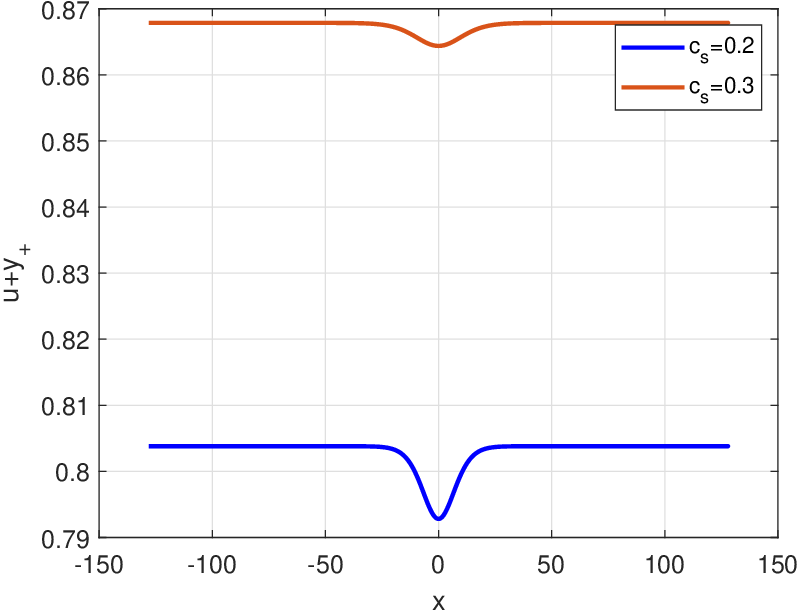}}
\subfigure[]
{\includegraphics[width=0.45\columnwidth]{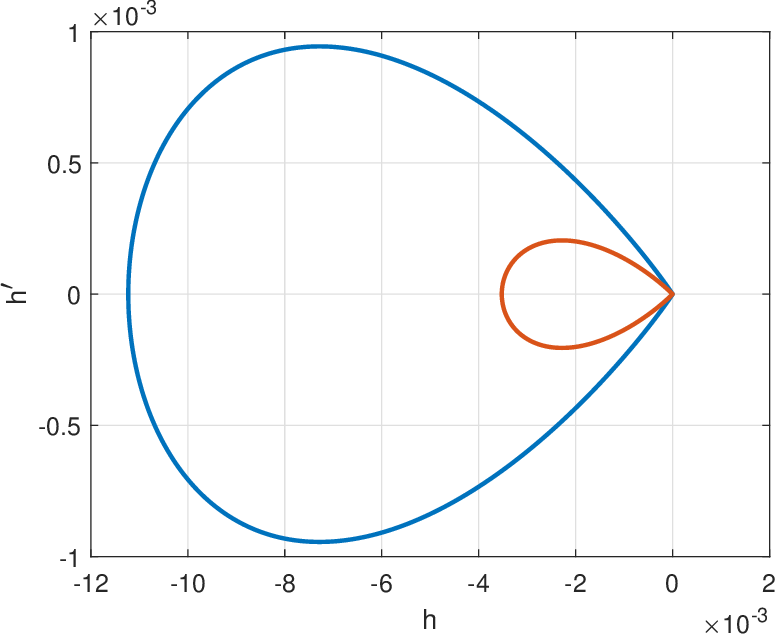}}
\subfigure[]
{\includegraphics[width=0.45\columnwidth]{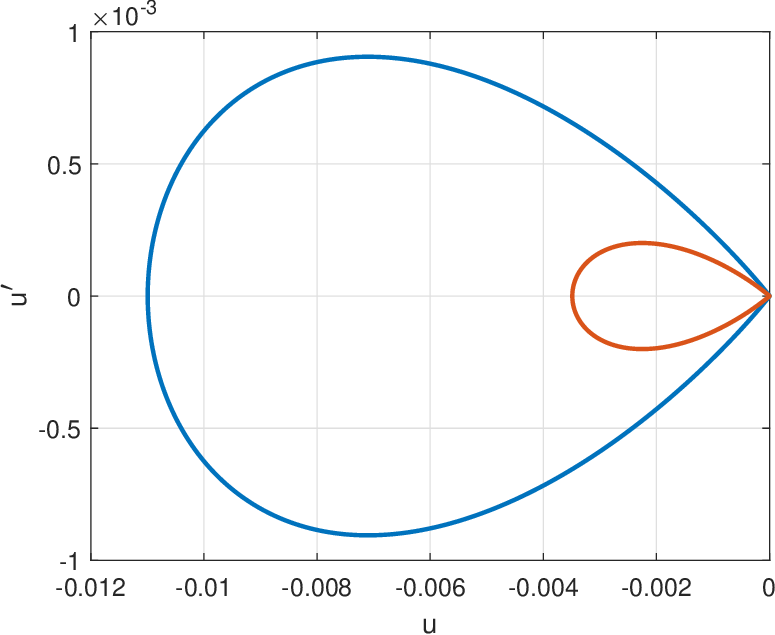}}
\caption{Approximations to traveling-wave solutions of (\ref{adg9}) for two values of $c_{s}$. Case $g=g_{-}(c_{s})+10^{-6}$. (a) $\widetilde{h}+y_{+}$ profile; (b) $\widetilde{u}+y_{+}$ profile; (c) $\widetilde{h}$ phase plot; (d) $\widetilde{u}$ phase plot.}
\label{ADG8a}
\end{figure}
\begin{figure}[htbp]
\centering
\centering
\subfigure[]
{\includegraphics[width=0.45\columnwidth]{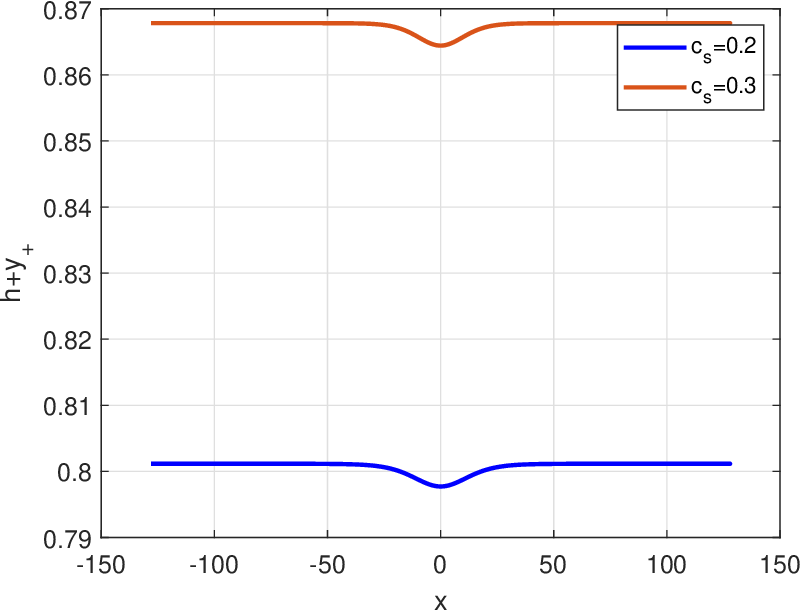}}
\subfigure[]
{\includegraphics[width=0.45\columnwidth]{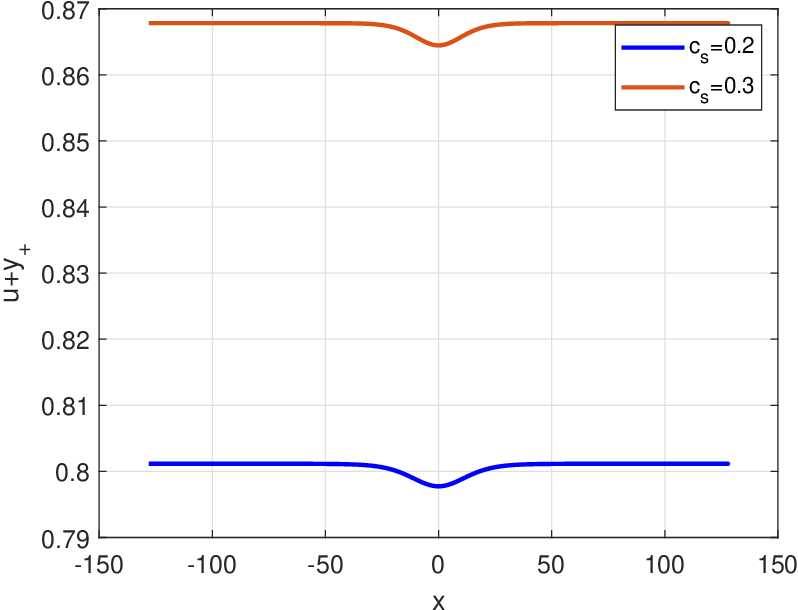}}
\subfigure[]
{\includegraphics[width=0.45\columnwidth]{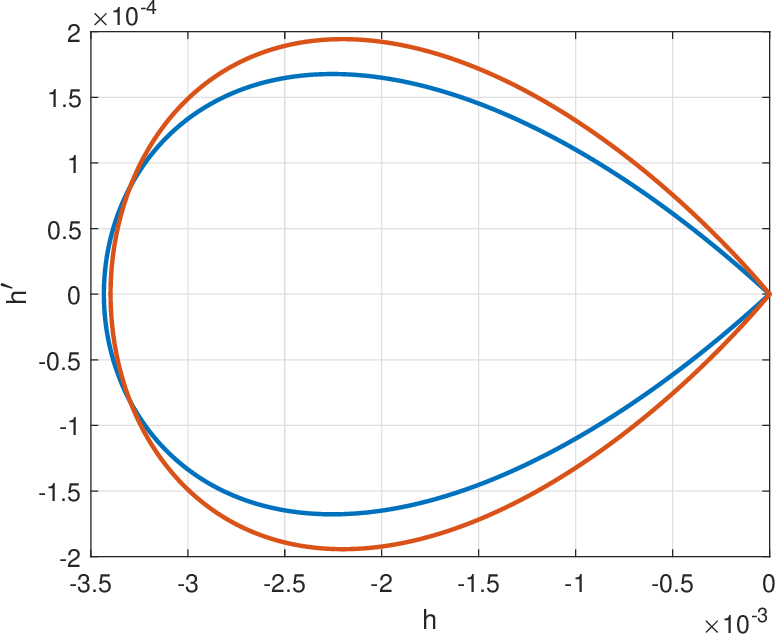}}
\subfigure[]
{\includegraphics[width=0.45\columnwidth]{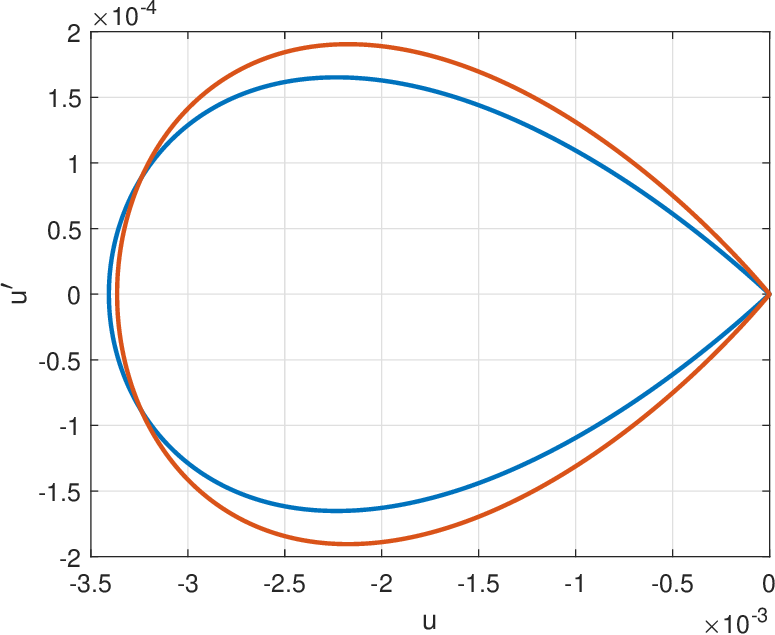}}
\caption{Approximations to traveling-wave solutions of (\ref{adg9}) for two values of $c_{s}$. Case $g=g_{1}(c_{s})-10^{-6}$. (a) $\widetilde{h}+y_{+}$ profile; (b) $\widetilde{u}+y_{+}$ profile; (c) $\widetilde{h}$ phase plot; (d) $\widetilde{u}$ phase plot.}
\label{ADG8b}
\end{figure}
\end{itemize}
The results are summarized in Table \ref{ADG_tav1}.

\begin{table}[ht]
\begin{center}
\begin{tabular}{|c|c|c|}
    \hline
$c_{s}>\frac{1}{2} (a_{+}>0)$&$(b,a)$ region&Type of wave\\\hline\hline
$g<g_{1}^{*}(c_{s})$&1R&NMTW (dep.)\\
$g_{1}^{*}(c_{s})\leq g<g_{-}(c_{s})$&1L&NMTW (dep.)\\
$g_{-}(c_{s})\leq g<g_{1}(c_{s})$&4&PTW\\\hline\hline
$c_{s}<\frac{1}{2} (a_{+}>0)$&$(b,a)$ region&Type of wave\\\hline\hline
$g<g_{-}(c_{s})$&1R&NMTW (dep.)\\
$g_{-}(c_{s})\leq g<g_{1}(c_{s})$&2&MTW (dep.)\\
    \hline\hline
\end{tabular}
\end{center}
\caption{Bifurcation study of homo-clinic orbits around the equilibrium $y_{+}$. NMTW (dep.): Nonmonotone traveling wave of depression; MTW: Monotone traveling wave of depression; PTW: Periodic traveling wave.}
\label{ADG_tav1}
\end{table}
\subsection{Dynamics around the equilibrium $y_{-}$}
\begin{figure}[htbp]
\centering
\centering
\subfigure[]
{\includegraphics[width=0.45\columnwidth]{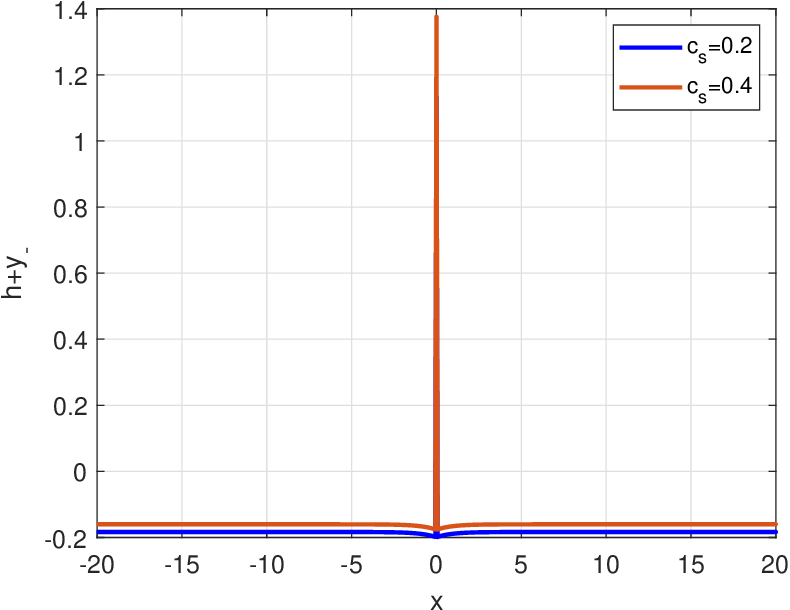}}
\subfigure[]
{\includegraphics[width=0.45\columnwidth]{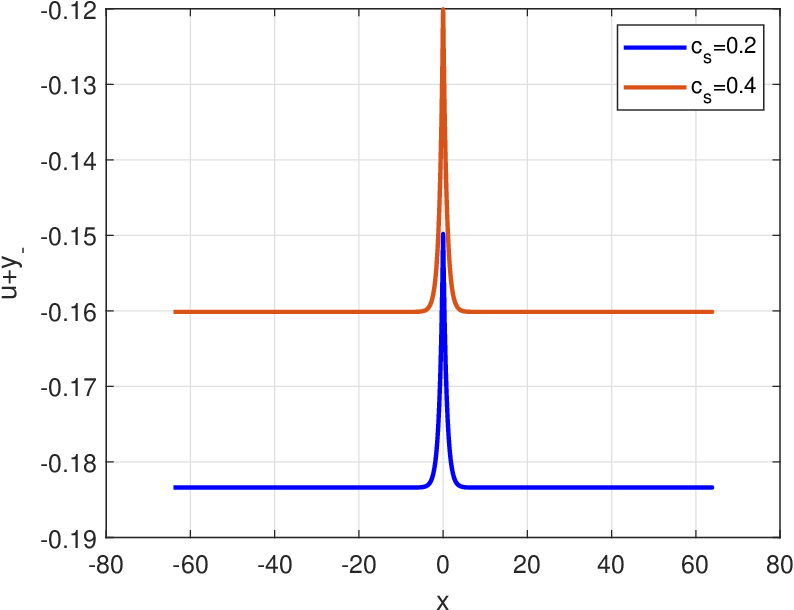}}
\caption{Approximations to traveling-wave solutions of (\ref{adg9}) for two values of $c_{s}$. Case $g=g_{+}(c_{s})-0.1$. (a) $\widetilde{h}+y_{-}$ profile; (b) $\widetilde{u}+y_{-}$ profile.}
\label{ADG12}
\end{figure}

\begin{figure}[htbp]
\centering
\centering
\subfigure[]
{\includegraphics[width=0.45\columnwidth]{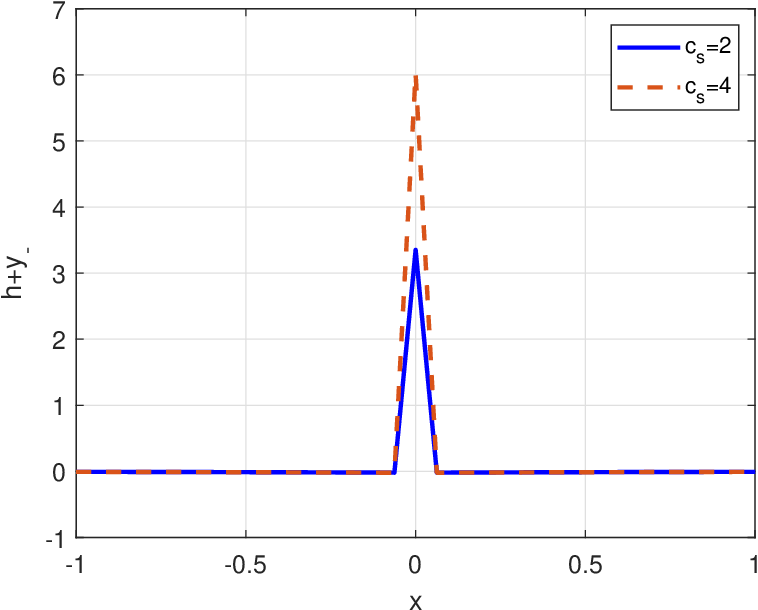}}
\subfigure[]
{\includegraphics[width=0.45\columnwidth]{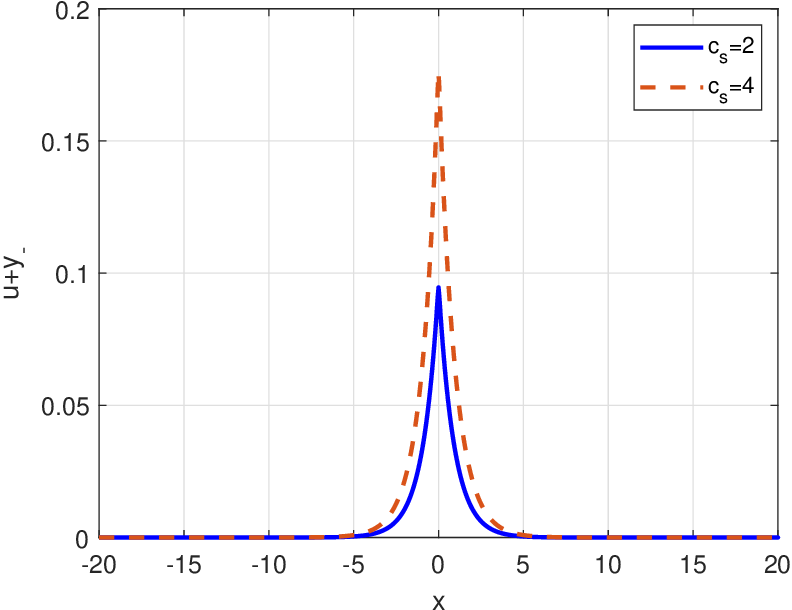}}
\caption{Approximations to traveling-wave solutions of (\ref{adg9}) for two values of $c_{s}$. Case $g=0$. (a) $\widetilde{h}+y_{-}$ profile; (b) $\widetilde{u}+y_{-}$ profile.}
\label{ADG13}
\end{figure}

\begin{figure}[htbp]
\centering
\centering
\subfigure[]
{\includegraphics[width=0.45\columnwidth]{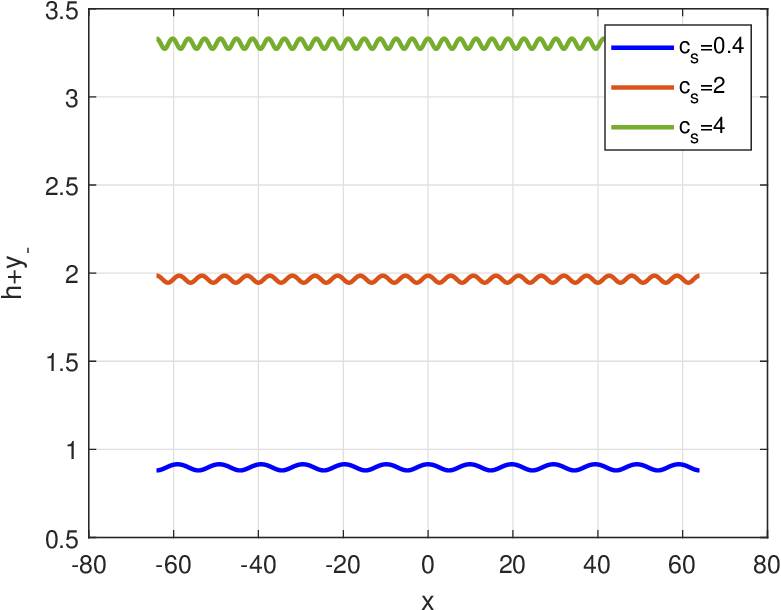}}
\subfigure[]
{\includegraphics[width=0.45\columnwidth]{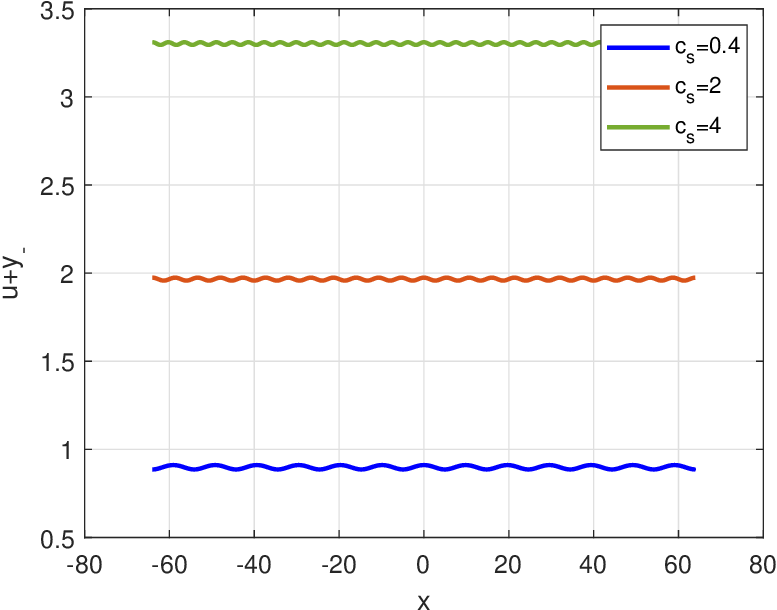}}
\caption{Approximations to traveling-wave solutions of (\ref{adg9}) for two values of $c_{s}$. Case $g=g_{1}(c_{s})-0.001$. (a) $\widetilde{h}+y_{-}$ profile; (b) $\widetilde{u}+y_{-}$ profile.}
\label{ADG14}
\end{figure}

A similar study for the equilibrium $y_{-}$ leads to the following results.
\begin{itemize}
\item[(i)] $c_{s}<\frac{1}{2}$. Here we may have four situations:
\begin{itemize}
\item[(i1)] $g<g_{+}(c_{s})$. Then $a_{-}, b_{-}>0$, and $\Delta_{-}<0$. Region 1, right.
\item[(i2)] $g_{+}(c_{s})\leq g<g_{1}^{**}(c_{s})$. Then $a_{-}, b_{-}>0$, and $\Delta_{-}>0$. Region 2.
\item[(i3)] $g_{1}^{**}(c_{s})\leq g<g_{1}^{*}(c_{s})$. Then $a_{-}, b_{-}<0$, and $\Delta_{-}>0$. Region 3, left.
\item[(i4)] $g_{1}^{*}(c_{s})\leq g<g_{1}(c_{s})$. Then $a_{-}<0, b_{-}>0$, and $\Delta_{-}>0$. Region 3, right.
\end{itemize}
\item[(ii)] $\frac{1}{2}<c_{s}<\frac{3}{}$. Here we may have four situations:
\begin{itemize}
\item[(ii1)] $g<g_{+}(c_{s})$. Then $a_{-}, b_{-}>0$, and $\Delta_{-}<0$. Region 1, right.
\item[(ii2)] $g_{+}(c_{s})\leq g<g_{1}^{**}(c_{s})$. Then $a_{-}, b_{-}>0$, and $\Delta_{-}>0$. Region 2.
\item[(ii3)] $g_{1}^{**}(c_{s})\leq g<g_{1}^{*}(c_{s})$. Then $a_{-}, b_{-}<0$, and $\Delta_{-}>0$. Region 3, left.
\item[(ii4)] $g_{1}^{*}(c_{s})\leq g<g_{1}(c_{s})$. Then $a_{-}<0, b_{-}<0$, and $\Delta_{-}>0$. Region 3, left.
\end{itemize}
\item[(iii)] $c_{s}>\frac{3}{}$. Here we may have four situations:
\begin{itemize}
\item[(iii1)] $g<g_{+}(c_{s})$. Then $a_{-}, b_{-}>0$, and $\Delta_{-}<0$. Region 1, right.
\item[(iii2)] $g_{+}(c_{s})\leq g<g_{1}^{*}(c_{s})$. Then $a_{-}, b_{-}>0$, and $\Delta_{-}>0$. Region 2.
\item[(iii3)] $g_{1}^{*}(c_{s})\leq g<g_{1}^{**}(c_{s})$. Then $a_{-}, b_{-}>0$, and $\Delta_{-}>0$. Region 2.
\item[(iii4)] $g_{1}^{*}(c_{s})\leq g<g_{1}(c_{s})$. Then $a_{-}<0, b_{-}<0$, and $\Delta_{-}>0$. Region 3, left.
\end{itemize}
\end{itemize}
The types of the generated profiles corresponding to the homo-clinic orbits with limit at $y_{-}$ are summarized in Table \ref{ADG_tav2}. Compared to Table \ref{ADG_tav1}, several differences can be mentioned:
\begin{itemize}
\item The types of the emerging traveling waves in the regions do not depend on the range (i), (ii) or (iii) of the speed $c_{s}$.
\item The non-periodic computed traveling waves are now of elevation.
\item In regions 1, right, and 2, classical traveling waves (of elevation) are  computed. The $h$ and $u$ profiles seem to have a non-monotone behaviour (more clearly observed in the first case) which seems to evolve to a monotone decay to the equilibrium $y_{-}$ as $g$ approaches $g_{+}(c_{s})$, cf. Figure \ref{ADG12}. (Thus, $(b_{-},a_{-})$ is closer to the curve $\mathbb{C}_{3}$ and eventually crosses from region 1, right, to region 2.) The profiles develop a peak at the point where the maximum is attained, which is more pronounced and larger for larger speeds, as observed in Figure \ref{ADG13}. Note in this case that $y_{-}=0$, and the computed traveling waves are solitary waves. A more rigorous proof of existence of such solitary waves may be considered by using classical approaches such as the concentration compactness theory of Lions, \cite{Lions}, or the positive operator theory, \cite{BenjaminBB1990}.
\item Crossing from region 2 to region 3 (left) generates a change from computed classical traveling waves to periodic traveling waves, see Figure \ref{ADG14}. PTW's are also computed in region 3, right.
\end{itemize}

\begin{table}[ht]
\begin{center}
\begin{tabular}{|c|c|c|}
    \hline
$c_{s}<\frac{1}{2}$&$(b,a)$ region&Type of wave\\\hline\hline
$g<g_{1}^{+}(c_{s})$&1R&TW (el.)\\
$g_{+}(c_{s})\leq g<g_{1}^{**}(c_{s})$&2&TW (el.)\\
$g_{1}^{**}(c_{s})\leq g<g_{1}^{*}(c_{s})$&3L&TW-PTW\\
$g_{1}^{*}(c_{s})\leq g<g_{1}(c_{s})$&3R&PTW\\\hline
\hline
$\frac{1}{2}<c_{s}<\frac{3}{}$&$(b,a)$ region&Type of wave\\\hline\hline
$g<g_{+}(c_{s})$&1R&TW (el.)\\\hline
$g_{+}(c_{s})\leq g<g_{1}^{**}(c_{s})$&2&TW (el.)\\
$g_{1}^{**}(c_{s})\leq g<g_{1}^{*}(c_{s})$&3L&TW-PTW\\
$g_{1}^{*}(c_{s})\leq g<g_{1}(c_{s})$&3L&TW-PTW\\\hline
    \hline
$c_{s}>\frac{3}{} (a_{+}>0)$&$(b,a)$ region&Type of wave\\\hline\hline
$g<g_{+}(c_{s})$&1R&TW (el.)\\\hline
$g_{+}(c_{s})\leq g<g_{1}^{*}(c_{s})$&2&TW (el.)\\
$g_{1}^{*}(c_{s})\leq g<g_{1}^{**}(c_{s})$&2&TW (el.)\\
$g_{1}^{**}(c_{s})\leq g<g_{1}(c_{s})$&3L&TW-PTW\\\hline
    \hline
\end{tabular}
\end{center}
\caption{Bifurcation study of homo-clinic orbits around the equilibrium $y_{-}$ TW (el.): Traveling waves of elevation; PTW: Periodic traveling waves.}
\label{ADG_tav2}
\end{table}

\section{Concluding remarks}
\label{sec4}

The present paper is concerned with the existence of traveling wave solutions of the asymptotic model (\ref{wave:uni}) for the evolution of a collision-free plasma in a magnetic field.  In Section \ref{sec:teo}, using Crandall-Rabinowitz theorem \cite{Crandall-Rabi} and the ideas in \cite{AGG}, we can rigorously prove the existence of traveling waves of small amplitude for (\ref{wave:uni}).

Furthermore, from a reformulation of the model which involves only local terms, the corresponding ode system for the profiles of the traveling wave solutions is derived and the existence of equilibria is discussed in terms of the speed of the wave and a constant of integration $g$, playing the r$\hat{o}$le of bifurcation parameter. The structure of solutions around two of these equilibria can be studied from that of the solutions of an equivalent first-order system (\ref{adg19}) around the origin as equilibrium. The reversible character of the system enables to analyze  the dynamics from the application of normal form theory for reversible vector fields, \cite{HaragusI}, in related references, \cite{Champ,ChampS,ChampT}. With this approach, and supported by an efficient numerical code to generate computed traveling wave profiles (described in Appendix \ref{appA}), the emergence of several traveling waves is revealed. The classification depends on the speed on the waves and the constant $g$ as bifurcation parameter. The types of traveling wave solutions include classical traveling waves of monotone and non-monotone behaviour, classical solitary waves (localized traveling waves) and periodic traveling waves. In addition, some experiments suggest the convergence of peakon-type waves, cf. \cite{ZhouT2010}.

\newpage

\appendix

\section{A procedure for the numerical generation of traveling waves}
\label{appA}
The numerical method to 
approximate solutions of (\ref{adg18}) is here formulated. The equation is first written in fixed-point form
\begin{equation}\label{unisw2}
\mathcal{L}\widetilde{u}=\mathcal{N}(\widetilde{u}),
\end{equation}
(we occasionally recover the tilde in the notation, for reasons that will become evident below) where $\mathcal{L}$ is a linear operator and $\mathcal{N}(\widetilde{u})$ gathers for the nonlinear terms involving $\widetilde{u}$ and its derivatives up to fourth order. Explicity
\begin{eqnarray*}
\mathcal{L}&=&\alpha_{\pm}\mathscr{J}^{2}+\frac{y_{\pm}}{2} \partial_{x}^{2}-\frac{1}{2}\mathscr{J},\nonumber\\
\mathcal{N}(\widetilde{u})&=&-\frac{1}{2}\left(\frac{3}{2}\mathscr{J}(\mathscr{J}{\widetilde{u}})^{2}+\partial_{x}(\widetilde{u}_{x}\mathscr{J}\widetilde{u})-\frac{1}{2}\mathscr{J}(\widetilde{u}_{x}^{2})\right).
\end{eqnarray*}
Equation (\ref{unisw2}) is iteratively solved from an approximation to (\ref{adg18}) on a long enough interval with periodic boundary conditions given by a Fourier collocation method. Let $N\geq 1$ be an even integer and let
\begin{equation}\label{sw3}
x_{j}=-l+j\Delta x,\quad j=0,\ldots,N-1, \Delta x=2l/N,
\end{equation}
be a uniform grid of collocation points on $(-l,l)$. We consider the finite dimensional space
$$S_{N}={\rm span}\{e^{ik\pi (x+l)/l}, k\in\mathbb{Z}, -\frac{N}{2}\leq k\leq \frac{N}{2}-1\},$$ and define the spectral Fourier collocation approximation to a solution $\widetilde{u}$ of the periodic problem associated to (\ref{adg18}) on $(-l,l)$ as $\widetilde{u}^{N}\in S_{N}$ satisfying (\ref{unisw2}) at the collocation points (\ref{sw3}). The approximation $\widetilde{u}^{N}$ is typically represented by the nodal values
\begin{equation}\label{sw4}
\widetilde{u}_{\Delta}=(\widetilde{u}_{0},\ldots,\widetilde{u}_{N-1})^{T},
\end{equation}
with $\widetilde{u}_{j}=\widetilde{u}^{N}(x_{j}), j=0,\ldots N-1$ as approximations to the nodal values of the solution $\widetilde{u}(x_{j})$. Thus the algebraic system satisfied by (\ref{sw4}) is of the form
\begin{equation}
\mathcal{L}_{\Delta}\widetilde{u}_{\Delta}=\mathcal{N}_{\Delta}(\widetilde{u}_{\Delta}),\label{sw5a}
\end{equation}
where
\begin{eqnarray*}
\mathcal{L}_{\Delta}&=&\alpha_{\pm}(I_{N}-D_{N}^{2})^{2}+\frac{y_{\pm}}{2}D_{N}^{2}-\frac{1}{2}(I_{N}-D_{N}^{2}),\nonumber\\
\mathcal{N}_{\Delta}(\widetilde{u}_{\Delta})&=&-\frac{1}{2}\left(\frac{3}{2}(I_{N}-D_{N}^{2})((I_{N}-D_{N}^{2}){\widetilde{u}_{\Delta}}).^{2}+D_{N}(D_{N}\widetilde{u}_{\Delta}.(I_{N}-D_{N}^{2})\widetilde{u}_{\Delta})\right.\\
&&\left.-\frac{1}{2}(I_{N}-D_{N}^{2})((D_{N}\widetilde{u}_{\Delta}).^{2})\right).
\end{eqnarray*}

where $I_{N}$ denotes the $N\times N$ identity matrix, $D_{N}$ stands for the pseudospectral differentiation matrix, \cite{CHQZ}, and the dots in the nonlinear part $\mathcal{N}_{\Delta}(\widetilde{u}_{\Delta})$ denotes Hadamard products of the corresponding vectors.
The fixed-point system (\ref{sw5a}) is implemented by writing its Fourier representation (in terms of the discrete Fourier components of the vectors (\ref{sw4}) via Discrete Fourier Transform) and solving iteratively the resulting algebraic system. A typical iterative resolution consists of combining the Petviashvili method with a vector extrapolation technique to accelerate the convergence. The details can be seen in e.~g. \cite{DDS2022} and references therein.

After obtaining approximations $\widetilde{u}_{\Delta}$ to the profile $\widetilde{u}$ solution of (\ref{adg18}), then we compute $u_{\Delta}=\widetilde{u}_{\Delta}+y$ to approximate (\ref{adg9c}) for both $y=y_{\pm}$ and, finally, we approximate the $h$ profile from the discrete version of (\ref{loc})
\begin{eqnarray*}
h_{\Delta}=(I_{N}-D_{N}^{2})u_{\Delta}.
\end{eqnarray*}
\end{document}